\documentclass[12pt]{article}
\usepackage{euscript,amsfonts,amssymb}
\author{A. Gorinov}
\title{Real cohomology groups
\date{}
of the space of nonsingular curves of degree 5 in
${\bf CP}^2$}
\begin{document}

\maketitle
\newcommand{\Z}{{\bf Z}}
\newcommand{\C}{{\bf C}}
\newcommand{\CP}{{\bf CP}}
\newcommand{\R}{{\bf R}}
\newcommand{\eu}{\EuScript}
\newcommand{\eps}{\varepsilon}
\newtheorem{theorem}{Theorem}
\newtheorem{lemma}{Lemma}
\newtheorem{Prop}{Proposition}
\newtheorem{definition}{Definition}
\newtheorem{corollary}{Corollary}

\section{Main result}

Let us denote by $\Pi_5$ the space of all homogeneous
polynomials $\C^3\to\C$ of degree 5
and by $P_5$ its subspace consisting of all nonsingular polynomials
(i.e. the polynomials, whose gradient is non-zero outside the origin).

\begin{theorem}
\label{main}
The Poincar{\'e}  polynomial of $P_5$ is equal to
$$(1+t)(1+t^3)(1+t^5).$$
\end{theorem}

A general method
of calculating the cohomology of the spaces of nonsingular
algebraic hypersurfaces
of given degree was described in \cite{vas1}.
In particular, the real cohomology groups of spaces of nonsingular
plane curves of degree $\le 4$ were calculated there. In the present
work, we apply a modification of the same method
to the case of nonsingular quintics in $\CP^2$.

Denote by $\Sigma_5$ the space $\Pi_5\setminus P_5$.
By the Alexander duality, the cohomology group
of $P_5$ is isomorphic to the Borel-Moore homology
(i.e. the homology of the complex of locally finite chains)
of $\Sigma_5$:
$$H^i(P_5,\R)=\bar H_{2D-1-i}(\Sigma_5,\R),$$
where $D=\dim_{\C}(\Pi_5)=21$ and $0<i<2D-1$. This reduction was used first
by V. I. Arnold in \cite{A}. To calculate the latter group
$\bar H_*(\Sigma,\R)$ we use a version of the spectral sequence
constructed in \cite{vas1}. It is described in the following theorem.

\begin{theorem}
\label{second}
The spectral sequence for the real Borel-Moore homology
of the space $\Sigma_5$ is defined by the following conditions:

1. Any its nontrivial term $E_{p,q}^1$ belongs to the
quadrilateral in the $(p,q)$-plane, defined by the conditions
$[1\le p\le 4, 28\le q \le 39]$.

2. In this quadrilateral all the nontrivial terms $E_{p,q}^1$
look as is shown in (\ref{spseq}).

3. The spectral sequence degenerates in this term, i.e.
$E^1 \equiv E^\infty.$
\end{theorem}

\begin{equation}
\label{spseq}
\begin{array}{ccccc}
39 & {\R}& & &\\
38 & & & &\\
37 & {\R} & &&\\
36 & & &&\\
35 & {\R}&{\R} &&\\
34 & & &&\\
33 & &{\R}& &\\
32 \\
31 & &{\R} & &\\
30 & & &&\\
29 & & &{\R}&\\
& 1 & 2 & 3 & \\
\end{array}
\end{equation}


The proof of Theorem \ref{second} will be given in Section \ref{quintics}

I wish to express my deep gratitude to
prof. V. A. Vassiliev for proposing the problem and stimulating
discussions and to A. A. Oblomkov and A. Yu. Inshakov, who explained me, how singular sets
of plane quintics look like.

\section{The method of conical resolutions}

Consider the following general situation: suppose $V$ is some
(finite-dimensional) vector space of (real or complex) functions on a manifold $\tilde M$ and $\Sigma\subset V$
is a discriminant formed by the functions that have singularities of a
given type. We want to calculate the Borel-Moore homology of $\Sigma$.
In order to do this we construct
a resolution, i.e., a space $\sigma$ and a proper map
$\pi:\sigma\to\Sigma$ such that the preimage of every point is
contractible. We are going to describe a construction of $\sigma$ via
configuration spaces. Our construction generalizes that from the article
\cite{vas1}.

Suppose that with every function $f\in\Sigma$ a compact
nonempty subset $K_f$ of some compact manifold $M$ is associated.
In the sequel we shall be interested in
the case when $\tilde M=\C^3\setminus\{0\}, V=\Pi_5, \Sigma=\Sigma_5$. In
this case it is natural to put $M$ equal to $\CP^2$ and $K_f$ equal to the
image of the set of singular points of $f$ under the evident map $\tilde
M\to\CP^2$.

In general, we suppose that the following conditions are satisfied: 1). If
$K_f\cap K_g\neq\varnothing$, then $f+g\in\Sigma$ and $K_f\cap K_g\subset
K_{f+g}$, 2). $K_{\lambda f}=K_f$ if $\lambda\neq 0$, 3). $K_0=M$, where
$0$ is the zero function.  Define for any $K\subset M$ the space
$L(K)\subset V$ as the space consisting of all $f$ such that $K\subset
K_f$. The conditions above imply that $L(K)$ is a vector space. For any
$x\in M$ we shall denote the space $L(\{x\})$ simply by $L(x)$.  We
suppose also that the dimension of $L(x)$ is the same for all $x\in M$,
say $d$, and that the map $x\mapsto L(x)$ from $M$ to $G_d(V)$ (the
corresponding Grassmann manifold) is continuous.

{\it Remark.} One may ask a natural question: if we are dealing with
functions on some manifold $\tilde M$, why should we introduce another manifold $M$?
The problem is that for our construction it will be convenient to
associate with a singular function a compact subset of a compact
manifold.  In the case when the manifold $\tilde M$ itself is compact, we
can assume, of course, that $M=\tilde M$, $K_f$ is the set of points where
$f$ has singularities of some given type.

By a configuration in a compact manifold $M$ we shall mean a compact
nonempty subset of $M$. Denote by $2^M$ the space of all
configurations in $M$. Suppose that the topology on $M$ is induced by a metric $\rho$.
We introduce the Hausdorff metric on
$2^M$ by the following rule:
$$\tilde\rho(K,L)=\mbox{max}_{x\in K}\rho(x,L)+\mbox{max}_{x\in
L}\rho(x,K).$$ Let us verify, for instance, that this definition satisfies
the triangle axiom:  Let $K, L, M\in X_i, \mbox{max}_{x\in
K}\rho(x,L)=\rho(x_0,y_0),$ $\mbox{max}_{y\in L}\rho(y,M)=\rho(y_1,z_1),$
$\mbox{max}\rho_{x\in K}\rho(x,M)=\rho(x_2,z_2), x_0, x_2\in K, y_0,
y_1\in L, z_1, z_2\in M$.  There exist $y'\in L, z'\in M$ such that
$\rho(x_2, y')\le\rho(x_0,y_0)$, $\rho(y',z')\le\rho(y_1,z_1)$. We have
$$\rho(x_0,y_0)+\rho(y_1,z_1)\ge\rho(x_2,y')+\rho(y',z')\ge\rho(x_2,z')\ge\rho(x_2,z_2).$$
Thus, we have proved that
$$\mbox{max}_{x\in K}\rho(x,L)+\mbox{max}_{y\in L}\rho(y,M)\ge\mbox{max}_{x\in K}\rho(x,M).$$ The inequality
$$\mbox{max}_{z\in M}\rho(z,L)+\mbox{max}_{y\in L}\rho(y,K)\ge\mbox{max}_{z\in M}\rho(z,K).$$ is obtained analogously.
Adding these two inequalities, we obtain the triangle inequality.

It is easy to check that if $M$ is compact, then the space $2^M$
provided with the metric $\tilde\rho$ is also compact. Let us denote by
$B(M,k)$ \label{bbb} the subspace of $2^M$ that consists of all configurations that contain exactly $k$ elements.
For any subspace
$A\subset 2^M$ we denote by $\bar A$ the closure of $A$ in $2^M$.
We have $\overline B(M,k)=\bigcup_{j\le k}B(M,j)$.

\begin{Prop}
Let $(K_j)$ be a Cauchy
sequence in $2^M$, and let $K$ be the set
consisting of the limits of all sequences
$(a_j)$ such that
every term $a_j\in K_j$. Then $K$ is nonempty and compact, and
$\lim_{j\to\infty}\tilde\rho(K_j,K)=0$.
\end{Prop}

$\diamondsuit$

\begin{Prop}
\label{inc}
Let $(K_i), (L_i)$ be two sequences in $2^M$. Suppose that there exist
$\lim_{i\to\infty}K_i,\lim_{i\to\infty}L_i$, and denote these limits by $K$
and $L$ respectively. Suppose also that $K_i\subset L_i$ for every $i$.
Then $K\subset L$.  \end{Prop}

$\diamondsuit$

Suppose that $X_1,\ldots,X_N$ are subspaces of $2^M$ that satisfy the
following conditions

\newcounter{ccc}
\label{conditions}
\begin{enumerate}
\item\label{first} For every $f\in\Sigma$ the set $K_f$
must belong to some $X_i$.

\item\label{sec} Suppose that $K\in X_i, L\in X_j, K\subsetneq L$. Then
$i<j$.

\item\label{three} Recall that $L(K)$ is the space of all
functions $f$ such that $K\subset K_f$. If we fix $i$, then the dimension of $L(K)$ must be the
same for all configurations $K\subset X_i$. We denote this dimension by
$d_i$.

\item\label{four} $X_i\cap X_j=\varnothing$ if $i\neq j$.

\item\label{five} Any $K\in\bar X_i\setminus X_i$ belongs to some $X_j$ with
$j<i$.

\item\label{six} For every $i$ the space $\eu T_i$ consisting of pairs $(x,K), x\in K, K\in
X_i$ is the total space of a locally trivial bundle over $X_i$ (the
projection $pr_i:\eu T_i\to X_i$ is evident). This
bundle will be called the {\it tautological bundle}\footnote{When,
for instance, $M=\CP^n$ and $X_i$ consists of projective subspaces of $M$
of one and the same dimension, this is just the projectivization of
the usual tautological bundle over some Grassmann manifold of $\C^{n+1}$.} over $X_i$.

\item\label{last}Note that any local trivialization of $\eu T_i$ has
the following form $(x,K')\mapsto (t(x,K'),K')$. Here $x$ is a point
in some $K\in X_i$, $K'$ belongs to some neighborhood $U\ni K$ in
$X_i$ and $t:K\times U\to M$ is a continuous map such that if we
fix $K'\in U$, we obtain a homeomorphism $t_{K'}:K\to K'$.
We suppose that for every $K\in X_i$ there exist a neighborhood $U\supset
K$ and a local trivialization of $\eu T_i$ over $U$ such that every map
$t_{K'}:K\to K'$ establishes a bijective correspondence between the
subsets of $K$ and $K'$ that belong to $\bigcup_{j\leq i}X_j$.
\addtocounter{ccc}{\value{enumi}}
\end{enumerate}

Under these assumptions we are going to construct a resolution $\sigma$ and a filtration on
it such that the $i$-th term of the filtration is the total space of a fiber
bundle over $X_i$.

Note that due to condition \ref{three} for any $i=1\ldots,N$ there exists an evident map $K\mapsto L(K)$ from $X_i$ to the
Grassmann manifold $G_{d_i}(V)$, which is continuous (since we
suppose that the map $x\mapsto L(x)$ is continuous).

{\it Remark.\label{6i7}} Note also that the rather strange-looking condition \ref{last} follows
immediately from condition \ref{six} in the following situation:
suppose $X_i$ consists of finite configurations, and for all $K,L$,
such that $K\in
X_i, L\subset K$ there is an index $j<i$ such that $L\in X_j$. In
this case any trivialization of $\eu T_i$ fits.

Consider the space $Y=\bigcup_{i=1}^N\bar
X_i=\bigcup_{i=1}^N X_i$. Denote by $X$ the space $Y^{*N}$ --- the $N$-th self-join of $Y$
\footnote{Recall that for any finite $CW$-complex $Y$ and any $k\ge 1$ the $k$-th self-join of $Y$
(denoted by $Y^{*k}$) can be defined as follows: take a generic
embedding $i:Y\to\R^{\Omega}$ for some
very large (but finite) $\Omega$ and define
$Y^{*k}$ as the union of all $(k-1)$-simplices with vertices in $i(Y)$
("generic" means that the intersection of any two
such simplices is their common face if there is any).}. Note that the
spaces $Y, Y^{*N}$ are compact.
Call a simplex $\triangle\subset X$ {\it coherent} if the configurations
corresponding to its vertices form an ascending chain. Note that then
its vertices belong to different $X_i$ (condition \ref{sec}). Among the
vertices of a coherent simplex there is one that contains as a subset
all the other.  Such vertex will be called the {\it main vertex} of
the corresponding simplex. Denote by $\Lambda$ the union of all coherent
simplices.  For any $K\in X_i$ denote by $\Lambda(K)$ the union of all
coherent simplices, whose main vertex is $K$. Note that the space
$\Lambda(K)$ is contractible.

Denote
by $\Phi_i$ the union $\bigcup_{j\leq i}\bigcup_{K\in X_j}\Lambda(K)$.
This is a filtration on $\Lambda:
\varnothing\subset\Phi_1\cdots\subset\Phi_N=\Lambda$.

For any simplex $\triangle\subset X$ denote by
$\stackrel{\circ}\triangle$ its interior, i.e., union of its points that
do not belong to the faces of lower dimension.  Note that for every $x\in
X$ there exists a unique simplex $\triangle$ such that
$x\in\stackrel{\circ}\triangle$.

\begin{Prop}
\label{predel}
Let $(x_i)$ be a sequence in $X$ such that
$\lim_{i\to\infty}x_i=x$. Suppose $x_i\in\triangle_i, x\in
\stackrel{\circ}\triangle$, where $\triangle_i$ are coherent
simplices, and suppose $K$ is a vertex of $\triangle$. Then there
exists a sequence $(K_i)$ such that $K_i$ is a vertex of
$\triangle_i$ and $\lim_{i\to\infty}K_i=K$.
\end{Prop}

$\diamondsuit$

\begin{Prop}
All spaces $\Lambda, \Lambda(K), \Phi_i$ are compact.
\end{Prop}

{\it Proof.} The proposition follows
immediately from Propositions \ref{predel} and \ref{inc}. $\diamondsuit$

For any $K\in X_i$ denote by $\partial\Lambda(K)$
the union $\bigcup_{\kappa}\Lambda(K)$
over all maximal $\kappa\in\cup_{j<i} X_j, \kappa\subsetneq K$. The
space $\Lambda(K)$ is the union of all segments that join points of
$\partial\Lambda(K)$ with $K$, and hence it is homeomorphic to the
cone over $\partial\Lambda(K)$.

Define the conical resolution $\sigma$ as the subspace of
$V\times\Lambda$ consisting of pairs $(f,x)$ such that
$f\in\Sigma, x\in\Lambda(K_f)$.
There exist
evident projections $\pi:\sigma\to\Sigma$ and $p:\sigma\to\Lambda$. We introduce a filtration on
$\sigma$ putting $F_i=p^{-1}(\Phi_i)$.
The map $\pi$ is proper, since the preimage of each compact set
$C\subset\Sigma$ is a closed subspace of $C\times\Lambda$, which is compact.

\begin{theorem}
\label{conical}
Suppose $X_1,\ldots,X_N$ are subspaces of $2^M$ that satisfy Conditions \ref{first}-\ref{last}. Then

\begin{enumerate}

\item $\pi$ induces an isomorphism of Borel-Moore homology groups.

\item Every space $F_i\setminus F_{i-1}$ is a complex vector bundle over
$\Phi_i\setminus\Phi_{i-1}$ of dimension $\dim(L(K)), K\in X_i$.

\item The space $\Phi_i\setminus\Phi_{i-1}$ is a fiber bundle over
$X_i$, the fiber being homeomorphic to
$\Lambda(K)\setminus\partial\Lambda(K)$.
\end{enumerate}
\end{theorem}

{\it Proof.}
The first statement of the theorem follows from the fact that $\pi$
is proper and $\pi^{-1}(f)=\Lambda(K_f)$, the latter space being contractible. To prove the second
statement, let us study the preimage of a point
$x\in\Phi_i\setminus\Phi_{i-1}$ under the map $p:F_i\setminus
F_{i-1}\to \Phi_i\setminus\Phi_{i-1}$. We claim that
$p^{-1}(x)=L(K)$ for some $K\in X_i$.

Recall that each point $x\in\Phi_i\setminus\Phi_{i-1}$ belongs to the
interior of some coherent simplex $\triangle$, whose main vertex lies
in $X_i$. Denote this vertex by $K$ and denote the map
$\Phi_i\setminus\Phi_{i-1}\ni x\mapsto K\in X_i$ by $f_i$.

Now suppose $f\in L(K)$, or, which is the same, $K_f\supset K$. This
implies $\Lambda(K)\subset\Lambda(K_f)$. So we have
$x\in\Lambda(K)\subset\Lambda(K_f)$ and $(f,x)\in\sigma, p(f,x)=x$,
hence $f\in p^{-1}(x)$.

Suppose now that $f\in p^{-1}(x)$. We have $(f,x)\in\sigma$, hence
$x\in\Lambda(K_f)$. This means that $x$ belongs to some coherent
simplex $\triangle'$, whose main vertex is $K_f$. But $x$ belongs to the
interior of $\triangle$, hence $\triangle$ is a face of $\triangle'$
and $K$ is a vertex of $\triangle'$. But $K_f$ is the main vertex of
$\triangle'$, hence $K\subset K_f$ and $f\in L(K)$.

So, we see that $p^{-1}(x)=L(K),K\in X_i$, and the second statement of
Theorem \ref{conical} follows immediately from the fact that the
dimension of $L(K)$ is the same for al $K\in X_i$ (condition \ref{three}).
In fact the bundle $p:F_i\setminus F_{i-1}\to\Phi_i\setminus\Phi_{i-1}$ is
the inverse image of the tautological bundle over $G_{d_i}(V)$ under the
map $x\stackrel{f_i}\mapsto K\mapsto L(K)$.

We shall show now that the map $f_i$ is a locally-trivial
fibration with the fiber $\Lambda(K)\setminus\partial\Lambda(K)$. This
will prove the third statement of Theorem \ref{conical}.

Denote by $\eu X_i$ the space consisting of all pairs $(L,K)$ such
that $L\subset K, L\in\bigcup_{j\le i}X_j$. Let $P:\eu X_i\to X_i$ be
the evident projection. We shall construct a trivialization of
$f_i:\Phi_i\setminus\Phi_{i-1}\to X_i$ from a particular trivialization of $\eu
T_i$ that exists due to condition \ref{last}.

For any $K\in X_i$ let $U$ and $t$ be the neighborhood of $K$ and the
trivialization of $\eu T_i$ over $U$ that exist due to condition \ref{last}.
For any $K,K'\in U, x\in K$, put $t_{K'}(x)=t(x,K')$.
Define a map $T:P^{-1}(K)\times U\to 2^M$ by the following rule: take
any $L\subset K, L\in\bigcup_{j\le i} X_j$ and put $T(L,K')$ equal to the
image of $L$ under $t_{K'}$. Due to condition \ref{last}
$T(\bullet,K')$ is a bijection between the set of $L\in 2^M$ such that
$L\subset K, L\in\bigcup_{j\le i}X_j$ and the set of $L'\in 2^M$ such that
$L'\subset K', L'\in\bigcup_{j\le i}X_j$.

\begin{Prop}
\begin{enumerate}

\item The map $T_1:P^{-1}(K)\times U\to 2^M\times U$ defined by the
formula $T_1(L,K')=(T(L,K'),K')$ maps $P^{-1}(K)\times U$
homeomorphically onto $P^{-1}(U)$.

\item For any $K'\in U$ we have $L\subset M\subset K$ iff $T(L,K')\subset
T(M,K')\subset T(K,K')=K'$.
\end{enumerate}
\end{Prop}

This implies that $\eu X_i$ is a locally trivial bundle over $X_i$.

{\it Proof.} The second statement is obvious. We have already seen that the
map $T_1$ is a bijection. The verification of the fact that $T_1$ and its
inverse are continuous is straightforward but rather boring. $\diamondsuit$

Now we can construct a trivialization of
$f_i:\Phi_i\setminus\Phi_{i-1}\to X_i$ over the same neighborhood
$U\ni K$
as above: take any $x\in f_i^{-1}(K)$
and any $K'\in U$. $x$ can be written in the form $x=\sum_j\alpha_j
L_j$, where all $\alpha_j>0, \sum\alpha_j=1, L_j\in\bigcup_{k\le
i}X_k$ and exactly one $L_j$ belongs to $X_i$.
Put $F(x,K')=\sum_k \alpha_k T(L_k,K')$. Again a straightforward
argument shows that $F$ is a homeomorphism $f_i^{-1}(K)\times U\to f_i^{-1}(U)$. Theorem
\ref{conical} is proved. $\diamondsuit$

Since every $\Lambda(K)$ is compact in the topology of
$\Lambda$, we
have $\bar{H}_*(\Lambda(K)\setminus\partial\Lambda(K))=
H_*(\Lambda(K),\partial\Lambda(K))=H_{*-1}(\partial\Lambda(K),pt)$.

Note conditions \ref{first}-\ref{last}
give us some freedom in choosing the $X_i$'s. For instance, we can
choose them so that the following condition would be
satisfied:

\begin{enumerate}
\addtocounter{enumi}{\value{ccc}}
\item\label{xxx}
if $K$ is a finite configuration from
$X_i$, then every subset $L\subset K$ is contained in some $X_j$ with
$j<i$. 
\end{enumerate}

We have the following two lemmas.

\begin{lemma}
\label{fin}
If Condition \ref{xxx} is satisfied, then the fiber
of the bundle $(\Phi_i \setminus \Phi_{i-1})\to X_i$ over a point
$K \in X_i$ is an open simplex, whose
vertices correspond to the points of the configuration $K$.
\end{lemma}
The proof is by induction on the number of points in $K\in X_i$. $\diamondsuit$

Note that in this situation the simplicial complex $\Lambda(K)$ is
(piecewise linearly) isomorphic
to the first barycentric subdivision of the simplex $\triangle$ spanned by the
vertices of $K$. The complex $\partial\Lambda(K)$ is isomorphic to the
first barycentric subdivision of the boundary of $\triangle$.

Moreover, denote by $\Lambda^{fin}$ the union of
$\Lambda(K)$ for all finite $K$.

\begin{lemma}
\label{fin1}
If Condition \ref{xxx} is satisfied, then there exists a map $C:\Lambda^{fin}\to
M^{*N}$ that maps every $K$ to an element of the
interior of the simplex $\triangle$
spanned by the points of $K$ that depends continuously on $K$.  This map is
a homeomorphism on its image, and it maps $\Lambda(K)$
(respectively, $\partial\Lambda(K)$) homeomorphically on $\triangle$ (respectibely, $\partial\triangle$).
\end{lemma}

$\diamondsuit$

It follows from Lemmas \ref{fin}, \ref{fin1} that for any $i$ such
that $X_i$ consists of finite configurations, the fiber bundle
$\Phi_i\setminus\Phi_{i-1}$ is isomorphic to the restriction to $X_i$ of the
evident bundle $M^{*k}\setminus M^{*(k-1)}\to B(M,k)$, where $k$ is
the number of points in a configuration from $X_i$. So we have
$$\bar H_*(\Phi_i\setminus\Phi_{i-1})=\bar H_{*-(k-1)}(X_i,\pm\R),$$
where $\pm\R$ is the coefficient system that will be defined a little later.
Recalling that all complex vector bundles are orientable and applying
statement 2 of Theorem \ref{conical}, we get
$$\bar H_*(F_i\setminus F_{i-1})=\bar
H_{*-2d_i}(\Phi_i\setminus\Phi_{i-1})=\bar
H_{*-2d_i-(k-1)}(X_i,\pm\R),$$
where $d_i=\dim_\C L(K),K\in X_i$.

\section{Some further preliminaries}

\begin{definition}
\label{def1}
For any space $X$ denote by $F(X,k)$ the space of all ordered $k$-ples from $X$, i.e. the space
$$X\times\cdots\times X\setminus \{(x_1,\ldots,x_k)|x_i=x_j\mbox{ for some }i\neq j\}.$$
\end{definition}

The space $B(X,k)$ defined on page \pageref{bbb} is the quotient of $F(X,k)$
by the evident action of the symmetric group $S_k$. We shall denote by $\tilde B(\CP^2,k)$ the
subspace of $B(\CP^2,k)$ consisting
of generic $k$-configurations (i.e. such configurations that contain no three points that belong
to a line, no 6 points that belong to a quadric etc.).
The spaces $\tilde F(\CP^2,k)$ are defined in a similar way.

\begin{definition}
\label{def2}
For any subspace $Y\subset B(X,k)$ of some configuration space denote by $\pm \R$
the "alternating" system, i.e.
the local system on $Y$ with the fiber $\R$ that changes it sign under the
action of any loop transposing exactly 2 points in a configuration from
$Y$.  \end{definition}

Throughout the text $P(X,{\eu L})$ stands for the Poincar{\'e} polynomial of $X$
"with coefficients
in ${\eu L}$", i.e. the polynomial $\sum_i a_it^i$, where
$a_i=\mbox{dim}(H^i(X,{\eu L}))$. In a similar way, we denote by $\bar P(X,{\eu L})$
the polynomial $\sum_i a_it^i$, where $a_i=\mbox{dim}(\bar{H}_i(X,{\eu
L}))$.

{\it We shall consider only homology and cohomology groups with real coefficients,
and the fibers of all local systems will be real vector spaces of finite dimension.}

The following statement will be frequently used in the sequel:

\begin{theorem}
\label{leray}
Let $p:E\to B$ be a locally trivial fiber bundle with
fiber $F$, and let ${\eu L}$ be a local system of groups or vector
spaces on $E$.  Then the Borel-Moore homology groups of $E$ with
coefficients in ${\eu L}$ can be obtained from the spectral sequence with
the term $E^2$ equal to $E^2_{p,q}=\bar H_p(B,{\bar\eu H}_q)$, where
${\bar\eu H}_q$ is the local system with the fiber $\bar H_q(F,{\eu L}|F)$
corresponding to the natural action of $\pi_1(B)$ on $\bar H_q(F,{\eu
L}|F)$.  \end{theorem} $\diamondsuit$ \medskip

This is a version of Leray theorem.
The fiber of $\bar \eu H_q$ equals $\bar H_q(F,\eu L|F)$ and a
loop $\gamma$ in $B$ defines a map $\tilde f:\eu L|F\to\eu L|F$ covering some
map $f:F\to F$. Recall the construction of $f$. To obtain $f$ consider
a family of curves $\gamma_x(t), x\in F$ that cover $\gamma$.
Then we can put $f(x)=\gamma_x(1), x\in F$. The map $\tilde f:\eu
L|F_*\to\eu L|F_*$ consists simply of the maps $\eu L_x\to\eu L_{f(x)}$ induced by
$\gamma_x$. For instance, suppose that $F$ is
connected and that $f$ preserves the distinguished point $*\in
F$. Then the restriction of $\tilde f$ to $\eu L_*$
is induced by the loop $t\mapsto\gamma_*(t)$.

The map $\tilde f$ induces for every $q$ a map $\tilde f_*:\bar H_q(F,\eu
L|F)\to \bar H_q(F,\eu L|F)$, which is exactly the map $\bar \eu H_q|_{\gamma(0)}\to\bar\eu H_q|_{\gamma(1)}$
induced by $\gamma$.

\begin{theorem}
\label{spseqmap}
Let $E_1
\to B_1, E_2\to B_2$ be two bundles and $\eu L_1, \eu L_2$ be
local coefficient systems of groups or vector spaces on $E_1, E_2$
respectively.  Let $f:E_1\to E_2$ be a proper map that covers some map
$g:B_1\to B_2$ (i.e., $f$ is a proper bundle map). Suppose $\tilde f:\eu
L_1\to\eu L_2$ is a map of coefficient systems that covers $f$.  Then the
map $\tilde f$ induces a map of the terms $E^2$ of spectral sequences of
Theorem \ref{leray}.  \end{theorem} $\diamondsuit$

The map of the spectral sequences induced by $\tilde f$ can be described explicitly in the
following way. Let $F^{(i)}_x$ be the fiber of $E_i$ over $x\in
B_i$, $\bar\eu H_q^{(i)}$ be the coefficient system on $B_i$
from Theorem \ref{leray}, $i=1,2$. Since $f$ is a bundle map, it
maps $F^{(1)}_x$ into $F^{(2)}_{g(x)}$, and $\tilde f$ maps $\eu
L|F^{(1)}_x$ into $\eu L|F^{(2)}_{g(x)}$. The latter map induces for any $x\in B_1$ a map
$\bar H_q(F^{(1)}_x,\eu L|F^{(1)}_x)\to\bar H_q(F^{(2)}_{g(x)},\eu
L|F^{(2)}_{g(x)})$, which can be considered as restriction to $\bar\eu H_q^{(1)}|_x$ of
a map $f'_q:\bar\eu
H_q^{(1)}\to\bar\eu H_q^{(2)}$ that covers $g$. The desired map of the
terms $E^2$ of
spectral sequences is just the map $\bar H_*(B_1,\bar\eu H_q^{(1)})\to
\bar H_*(B_2,\bar\eu H_q^{(2)})$ induced by $f'$.

Theorem \ref{leray} has the following corollary:
\begin{corollary}
\label{cor1}
Let $N,\tilde N$ be manifolds, let $p:\tilde N\to N$ be a finite-sheeted
covering, and let ${\eu L}$ be a local system of groups on $\tilde N$. Then
$H^*(\tilde N,{\eu L})=H^*(N,p(\eu L))$, where $p(\eu L)$ denotes the
direct image of the system $\eu L$.  \end{corollary} $\diamondsuit$


If $\eu L$ is the constant local system with fiber $\R$,
and the map $\tilde N\to N$ is the quotient by a free
action of a finite group $G$, then the representation of
$\pi_1(N)$ on a fiber of $p(\eu L)$ is isomorphic to the natural action of
$\pi_1(N)$ on the group algebra of $\pi_1(N)/p(\pi_1(\tilde N))\cong G$.
In particular if $\tilde N$ is simply connected, we get the regular representation of $G$
(the group acts on its group algebra by left shifts).

Recall that any irreducible real or complex representation of a finite
group $G$ is included into the regular representation. If the
representation is complex, it is obvious. In the real case it follows from
Schur's Lemma and from the well-known fact that if $R:G\to GL(V), S:G\to
GL(W)$ are real representations of $G$ and $R_{\C}, S_{\C}$ are their
complexifications, then $\mbox{dim}_{\R}(\mbox{Hom}(R,S))=
\mbox{dim}_{\C}(\mbox{Hom}(R_{\C},S_{\C}))$, where $\mbox{Hom}(R,S)$ is the space of
representation homomorphisms between $R$ and $S$ (i.e. such operators $f:V\to W$ that
$f(R(g)x)=S(g)f(x)$ for any $g\in G, x\in V$).

The homological analogues of Theorems \ref{leray} and \ref{spseqmap} are
obtained by omitting all bars over $H$-es and $\eu H$-es.
The cohomological versions of these theorems are obtained as follows:
in Theorem \ref{leray} the
action of a loop $\gamma$ is the inverse of the cohomology map induced by $\tilde
f:\eu L_x\to\eu L_{f(x)}$, and in Theorem \ref{spseqmap} we have
to suppose that $\eu L_1$ is the inverse image of $\eu L_2$, i.e. that the
restriction of $\tilde f$ over each point is bijective.

Neither in homological, nor in cohomological analogue of Theorem \ref{spseqmap} we have to require that
$f$ is proper.

We shall also need the following version of Poincar{\'e} duality theorem:

\begin{theorem}
\label{poincare}
Let $M$ be a manifold of dimension $n$, and let $\eu L$ be a local system on $M$, whose fiber is a
real or complex vector space. Then we have
$$H^*(M,\eu L\otimes Or(M))\cong\bar H_{n-*}(M,\eu L),$$
where $Or(M)$ is
the orienting sheaf of $M$.
\end{theorem}
$\diamondsuit$

The following lemma allows us to calculate the real cohomology groups of the quotient of a semisimple connected
Lie group by a finite subgroup
with coefficients in arbitrary local systems. We shall use this lemma several times.

\begin{lemma}
\label{lem3}
Let $G$ be a connected and simply connected Lie group that is either compact or complex
semisimple.
Let $G_1$ be a finite subgroup of $G$, and let $G/G_1$ be the coset of left
classes by $G_1$.
The cohomology group
$H^*(G/G_1,\eu L)$ is trivial if
the system $\eu L$ is nontrivial
and the action of
$\pi_1(G/G_1)\cong G_1$ on the fiber of ${\eu L}$ is irreducible.
\end{lemma}

(Here and below the "trivial system" means the constant sheaf with the fiber $\R$. In the sequel we shall denote this system simply by $\R$.)

{\it Proof.} If we apply Corollary \ref{cor1} to the covering
$G\to G/G_1$, we obtain

$$H^*(G,\R)=\oplus_i H^*(G/G_1,\eu L_i), $$
where $\eu L_i$ are coefficient systems corresponding to
irreducible representations of $G_1$ that are included into the regular real
representation.
Recall that in fact all real irreducible representations of a finite group are included into
the regular real representation.

But $H^*(G,\R)=H^*(G/G_1,\R)$, since every cohomology class of $G$ can be represented
by an invariant form,
therefore all the groups $H^*(G/G_1,{\eu L}_i)$ for $\eu L_i\neq\R$
are zero.
$\diamondsuit$

After Theorem \ref{leray} we gave an explicit construction of local systems that appear in the Leray spectral sequence
of a locally-trivial fibration. However, in many interesting cases we have a map that is ``almost'' a fibration, say the quotient
of a smooth manifold by an almost free action of a compact group etc., and we would like to know, how the Leray sequence
(which is defined for any continuous map) looks like in this case. It turns out that in the case of a quotient map the sheaves that
occur in the Leray sequence can be explicitly described (at least for some actions of some groups).

Let us fix a smooth action of a Lie group $G$ on a manifold $M$.
A submanifold $S\subset M$ is called a slice at $x\in S$ for the action of $G$ iff
$GS$ is open in $M$, and
there is a $G$-equivariant map $GS\to G/Stab(x)$ such  that the preimage of $Stab(x)$ under this map is
$S$. Here $GS$ is the union of the orbits of all points of $S$, and $Stab(x)$ is the stabilizer of $x$.

If $G$ is compact, a slice exists for any action at every point $x\in M$: provide $M$ with a $G$-invariant Riemannian metric
and put $S$ equal to the exponential of an $\varepsilon$-neighborhood of zero in the orthogonal
complement of $T_x(Gx)$ for any sufficiently small $\varepsilon$. (Here $Gx$ is the orbit of $x$.)

\begin{theorem}
\label{groupactions}
Suppose that $G$ is compact or complex semisimple. Suppose that $G$ is connected, simply connected
and for any $x\in M$ the group $Stab(x)$ is finite. Let $\eu L$ be a local system on $M$.
If there exists a slice for the action of $G$ at every point $x\in M$, then for any $i$ the sheaf $\eu H^i_{\eu L}$ on $M/G$
generated by the presheaf $U\mapsto H^i(p^{-1}(U),\eu L)$ is isomorphic to $p(\eu L)\otimes H^i(G,\R)$, where $p:M\to M/G$ is the
natural projection, $p(\eu L)$ is the direct image of $\eu L$, and $H^i(G,\R)$ is the constant sheaf with fiber $H^i(G,\R)$.
Moreover if $\eu L=\R$, then $p(\eu L)=\R$.
\end{theorem}

Note that the sheaf $\eu H^0_{\eu L}$ is canonically isomorphic to $p(\eu L)$ for any $\eu L$.

{\it Proof.}
Let us first consider the case $\eu L=\R$. We have to show that for any $i$ the sheaf $\eu H^i_\R$ is constant with fiber
isomorphic to $H^i(G,\R)$. Let $S$ be a slice for the action of $G$ at $x\in M$. It is easy to show that $S$ is invariant with respect to
$Stab(x)$ and $GS$ is homeomorphic to $G\times_{Stab(x)}S$, which is the quotient of $G\times S$ by the following action of
$Stab(x):g(g_1,x_1)=(g_1g^{-1},gx_1)$ for any $g\in Stab(x), g_1\in G, x_1\in S$.
This imply easily that 
for any $x'\in M/G$ open sets $U\ni x'$ such that
$p^{-1}(U)$ contracts to $p^{-1}(x')$ form a local basis at $x'$. Hence the canonical map $\rho_{x'}:\eu
H^i_\R(x')\to H^i(p^{-1}(x'),\R)$ is an isomorphism, where $\eu H^i_\R(x)$ is the fiber of $\eu H^i_\R$ over $x$.

Note that for any $x\in M$ the action map $\tau_x:G\to Gx, \tau_x(g)=gx$ induces
an isomorphism of real cohomology groups (the existence of a slice
at $x$ implies that $Gx$ is homeomorphic to $G/Stab(x)$, and, since $G$ is complex semisimple or compact,
and $Stab(x)$ is finite, the real cohomology map induced by $G\to G/Stab(x)$ is an isomorphism). Let $\sigma:M/G\to M$ be any map,
such that $p\circ\sigma=Id_{M/G}$.
Now for any $i, y\in H^i(G,\R)$ define the section $s_y$ of $\eu H^i_\R$ as follows: $s_y(x')=\rho^{-1}_{x'}\circ (\tau_{\sigma(x')}^*)^{-1}(y)$.

Let $U$ be a neighborhood of $x'$ such that $p^{-1}(U)$ contracts to $p^{-1}(x')$. It is easy to check that for any $y \in H^i(G,\R)$ there
is a $y'\in H^i(p^{-1}(U))$ such that $s_y$ coincides on $U$ with the canonical section of $\eu H^i_\R$ over $U$ defined by $y'$, so all maps
$x'\mapsto s_y(x')$ are indeed sections. It follows immediately from the definition that the section $s_y$ is nowhere zero if $y\neq 0$.
Putting $y$ equal to elements of some basis in $H^i(G,\R)$, we obtain $\dim(H^i(G,\R))$ everywhere linearly independent sections of $\eu H^i_\R$, whose
values span $\eu H^i_\R(x')$ for any $x'$. Hence the map $(x',y)\mapsto f_{x'}(y)$ establishes an isomorphism between $\eu H^i_\R$ and $H^i(G,\R)$.
The theorem is proved in the case, when $\eu L$ is trivial.

Now suppose that the system $\eu L$ is arbitrary. Note that for any open subset $U$ of $M/G$ there is a natural map (the $\smile$-product)
$$H^0(p^{-1}(U),\eu L)\otimes H^i(p^{-1}(U),\R)\to H^i(p^{-1}(U),\eu L\otimes\R)\cong H^i(p^{-1}(U), \eu L)$$
This gives us a map
\begin{equation}
\label{smile}
\eu H^0_{\eu L}\otimes \eu H^i_\R\to\eu H^i_{\eu L}.
\end{equation}

Due to the existence of a slice at each point of $M$ the groups
$\eu H^0_{\eu L}(x'),\eu H^i_\R(x'), \eu H^i_{\eu L}(x')$ are canonically isomorphic to $H^0(p^{-1}(x'),\eu L), H^i(p^{-1}(x'),\R),
H^i(p^{-1}(x'), \eu L)$ respectively (for any $x'\in M/G$). Under this identification the restriction of the map (\ref{smile}) to the fiber over
a point $x'\in M/G$ is the cup product map
\begin{equation}
\label{smile1}
H^0(p^{-1}(x'),\eu L)\otimes H^i(p^{-1}(x'),\R)\to H^i(p^{-1}(x'), \eu L).
\end{equation}
Now take some $x\in p^{-1}(x')$.
The action of $\pi_1(p^{-1}(x'))\cong Stab(x)$ in the fiber of $\eu L|p^{-1}(x')$ splits into a sum of irreducible representations, hence $\eu L|p^{-1}(x')$
can be decomposed into a sum $\eu L=\oplus\eu L_j$ such that the action of $\pi_1(p^{-1}(x'))$ on a fiber of each $\eu L_j$ is irreducible. Note that the map
(\ref{smile1}) respects the decomposition of $\eu L$ into a direct sum. Applying Lemma \ref{lem3} to each $\eu L_j$, we see that (\ref{smile}) is an
isomorphism. We have already shown that $\eu H^i_\R$ is constant with fiber isomorphic to $H^i(G,\R)$. The theorem is proved.$\diamondsuit$

The following 3 lemmas will be frequently used in our calculations (see
\cite{vas1}, \cite{vas2} for a proof):

\begin{lemma}
The group $\bar H_*(B(\C^n,k),\pm\R)$ is trivial for any $k\ge 2,
n\ge1$.
\end{lemma}

$\diamondsuit$

\begin{lemma}
The group $\bar H_*(B(\CP^n,k),\pm\R)$
for $n\ge 1$ is isomorphic to $H_{*-k(k-1)} (G_k(\C^{n+1}), \R)$,
where $G_k(\C^{n+1})$ is the Grassmann manifold of
$k$-dimensional subspaces in $\C^{n+1}$.
\end{lemma}
In particular, the group  $\bar H_*(B(\CP^n,k),\pm\R)$ is trivial if $k>n+1$.

$\diamondsuit$

\begin{lemma}
If $k\ge 2$, then the group $H_*((S^2)^{*k},\R)$
is trivial in all positive dimensions, where $(S^2)^{*k}$ is the
$k$-th self-join of $S^2$.
\end{lemma}

$\diamondsuit$

Consider the space $\C\setminus\{\pm 1\}$ and the coefficient system
$\eu L$ on it that changes its sign under any loop based at $0$ that passes once around $1$ or $-1$.
Let $f$ be the map $z\mapsto -z$ and let
$\tilde f:\eu L\to\eu L$ be the map that covers $f$ and is identical over
$0$.

\begin{Prop}
\label{C}
The map $\tilde f$ multiplies by $-1$ the groups $H_1(\C\setminus\{\pm
1\},\eu L)$, $H^1(\C\setminus\{\pm 1\},\eu L)$, $\bar H_1(\C\setminus\{\pm
1\},\eu L)$.
\end{Prop}

$\diamondsuit$

Consider $B(\C^*,2)$ --- the space of
couples of points in $\C\setminus\{0\}$. It is a fiber bundle over
$\C^*$, the projection $p:B(\C^*,2)\to\C^*$ being just the multiplication. The fiber is
homeomorphic to $\C\setminus\{\pm 1\}$, and the action of the generator of
$\pi_1(\C^*)$ is $z\mapsto -z$. The fiber $p^{-1}(1)$ contracts to the caracter ``8''.
Denote by $b, c$ the loops in $p^{-1}(1)$ corresponding to the circles of the ``8'', and denote
by $a$ the loop $t\mapsto\{ ie^{\pi it},-ie^{\pi it}\}$.

Consider the following three local
systems on $B(\C^*,2)$ (the fiber of each of them is $\R$):
\begin{enumerate}
\item $\eu A_1$
changes its sign under $a$
and does not change its sign under $b$ and $c$.

\item $\eu A_2$ changes its sign  under $b$ and $c$ and does
not change its sign under $a$.

\item $\eu A_3$ changes its sign under all loops $a,b,c$.
\end{enumerate}
Let $f$ be the map $B(\C^*,2)\to B(\C^*,2)$ induced by the map
$z\mapsto 1/z$ and let $f^i:\eu A_i\to \eu A_i, i=1,2,3$ be the
map that covers $f$ and is identical over the fiber $p^{-1}(1)$.

\begin{Prop}
\label{B}
\begin{enumerate}
\item $\bar P(B(\C^*,2),\eu A_1)=\bar P(B(\C^*,2),\eu A_3)=t^2(1+t),$
$\bar P(B(\C^*,2),\eu A_2)=0$.
\item The map $f^i_*$ preserves $\bar H_3(B(\C^*,2),\eu A_i)$ and
multiplies $\bar H_2(B(\C^*,2),\eu A_i)$ by $-1$ (i=1,3).
\end{enumerate}
\end{Prop}

Note that the lifting of the generator of $\pi_1(\C^*)$ looks like:
$\gamma_{\{a,b\}}(t)=\{ae^{\pi it},be^{\pi it}\}$, where $\{a,b\}\in p^{-1}(1)$.

This statement follows immediately from Theorems \ref{leray},
\ref{spseqmap}. $\diamondsuit$

\section{The spectral sequence for plane quintics}
\label{quintics}
Now we put $V=\Pi_5, \Sigma=\Sigma_5$.

\begin{Prop}
\label{class}
The configuration spaces $X_1,\ldots,X_{\ref{cp2}}$ that consist of the
following configurations satisfy Conditions \ref{first}-\ref{last} and Condition \ref{xxx}. The number indicated in square brackets equals
the dimension of $L(K)$ for $K$ lying in the corresponding $X_i$.

\begin{enumerate}
\item One point in $\CP^2$ [18].

\item 2 points in $\CP^2$ [15].

\item 3 points points in $\CP^2$ [12].

\item\label{4ptslin} 4 points on a line [11].

\item\label{5ptslin} 5 points on a line [10].

\item 6 points on a line [10].

\item 7 points on a line [10].

\item\label{8lin} 8 points on a line [10].

\item 9 points on a line [10].

\item\label{10lin} 10 points on a line [10].

\item\label{line} A line in $\CP^2$ [10].

\item\label{4nlin} 4 points in $\CP^2$ not on a line [9]. (Any three of them may
belong to a line though.)

\item\label{4+1} 4 points on a line $+$ a point not belonging to the line [8].

\item\label{5+1} 5 points on a line $+$ one point not belonging to the line [7].

\item 6 points on a line $+$ one point not belonging to the line [7].

\item\label{7+1} 7 points on a line $+$ one point not belonging to the line [7].

\item\label{line+1} A line in $\CP^2+$ a point not belonging to the line [7].

\item\label{5gen} 5 points in $\CP^2$ such that there is no line containing 4 of them [6].

\item\label{4+2} 4 points on a line $+$ 2 points not belonging to the line [5].

\item\label{5+2} 5 points on a line $+$ 2 points not belonging to the line [4].

\item\label{6+2} 6 points on a line $+$ 2 points not belonging to the line [4].

\item\label{line+2} A line in $\CP^2+$ 2 points not belonging to it [4].

\item \label{3+3} 3 points on each of two intersecting lines such that none of the
points coincides with the point of intersection [4].

\item\label{6q} 6 points on a
nondegenerate quadric [4].

\item\label{3+3+int} A configuration of type \ref{3+3} $+$ the point of intersection
of the lines [4].

\item\label{6gen} 6 points not belonging to a (possibly degenerate) quadric such that there is no
line containing 4 of those points [3].

\item\label{4+3lin} 4 points on a line $+$ 3 points on another line such that none
of the 7 points coincides with the point of intersection of the lines [3].

\item\label{5+3lin} 5 points on a line $l_1$ $+$ 3 points on $l_2\setminus l_1$, where $l_2$ is a line $\neq l_1$ [3].

\item\label{line+3line} A line $+$ 3 points of some other line, none of which coincides
with the point of intersection of the lines [3].

\item\label{4+4} 4 points on a line $+$ 4 points on some other line such that none of the 8
points coincides with the point of intersection of the lines [3].

\item\label{2lines} A couple of lines in $\CP^2$ [3].

\item\label{7q} 7 points on a nondegenerate quadric [3].

\item\label{quadric} A nondegenerate quadric [3].

\item\label{4+3} 4 points on a line $+$ 3 generic points not belonging to the line [2].

\item\label{3+3+1} A configuration of type \ref{3+3} $+$ a point not belonging to
the union of the lines [1].

\item\label{6q+1} 6 points on a nondegenerate quadric $+$ a point
not belonging to the quadric [1].

\item\label{3+3+1+int} A configuration of type \ref{3+3+1}  $+$ the point of intersection of the lines [1].

\item\label{8} 4 generic points $\{A, B, C, D\}$
in ${\CP^2}+$ 4 points of intersection of a line $l$ not passing through $A, B, C, D$
and two (possibly degenerate) quadrics passing through $A, B, C, D$ and not tangential to $l$ [1].

\item\label{9} 3 generic points $\{A, B, C\}$ in
${\CP^2}+$ 6 points of intersection of 3 lines $AB, BC, AC$ and
a (possibly degenerate) quadric not passing through $A,B,C,$ and
not tangential to the lines $AB, BC, AC$ [1].

\item\label{10} 10 points of intersection of 5 generic lines in $\CP^2$ [1].

\item\label{line+3} A line in $\CP^2+$ 3 generic points not belonging to the line [1].

\item\label{cp2} The whole $\CP^2$ [0].
\end{enumerate}

\end{Prop}
$\diamondsuit$
\medskip

{\it Proof. } We shall say that points $x_1,\ldots, x_k\in\CP^2$ are generic, iff among these points there
are no 3 points that are on the same line, no 6 points on a quadric, no 10 points on a cubic, etc.

Condition \ref{first} follows from the following observations: 1). the singular set
of a curve defined by a product of two polynomials is the union of the singular sets of the
curves defined by those polynomials and the intersection points of the curves; 2). the singular set
of an irreducible curve of degree 5 consists of 1 to 6 generic points
3). all possible singular sets of curves of
degree $\leq 4$ are described in \cite{vas1}. Note that some spaces $X_i$ contain
(or consist of) configurations that are not sets of singular points of any curve of degree 5. We introduce them to
make sure that Conditions \ref{five} and \ref{xxx} are satisfied.

The verification of Conditions \ref{sec}, \ref{four} and \ref{xxx} is straitforward.

Condition \ref{three} can be deduced from the following lemma.

\begin{lemma}
\label{lll}
\begin{enumerate}
\item\label{xyz1} Let $x_1,\ldots,x_k, k\leq 6$ be generic points of $\CP^2$.
Then the complex dimension of the space $L(\{x_1,\ldots,x_k\})$ (which consists of polynomials of degree 5 that have
singularities in all points $x_1,\ldots,x_k$ and maybe, elsewhere) is $21-3k$.

\item\label{xyz2} Let $l_1,l_2$ be two distinct lines in $\CP^2$. Suppose that  $x_1,x_2,x_3\in l_1\setminus l_2, y_1,y_2,y_3\in
l_2\setminus l_1, A\not\in l_1\cup l_2$. Then there exists exactly one cubic passing through all the points $x_i, y_j, i,j=1,2,3$ and
having a singularity at $A$.

\item\label{xyz3} Let $Q$ be a nondegenerate quadric in $\CP^2$, and suppose that $x_1,\ldots,x_6\in Q, A\not\in Q$. Then there exists
exactly one cubic passing through $x_1,\ldots,x_6$ and having a singularity at $A$.

\item\label{xyz4} If a curve of degree 5 has 3 singular points on a line, then it contains the line. If a curve of degree 5 has 6
singular points on a nondegenerate quadric, then it contains the quadric.

\item\label{xyz5} Consider a point $A\in\CP^2$. For any $d$ define $L^x_d(A)$ (respectively, $L^y_d(A), L^z_d(A), M_d(A)$) as the linear
space of
homogeneous polynomials $f$ of degree $d$ such that $\frac{\partial f}{\partial x}=0$ (respectively, $\frac{\partial f}{\partial y}=0$,
$\frac{\partial f}{\partial z}=0$, $f=0$) at every point of the preimage of $A$ under the natural map $\C^3\setminus\{0\}\to\CP^2$. Suppose
that $l$ is a line in $\CP^2, x_1,x_2,x_3,x_4\in l, y_1,y_2,y_3\not\in l$, and suppose that $y_1,y_2,y_3$ are not on a line. Then
$$dim_\C((\cap_{i=1}^4 M_4(x_i))\cap (\cap_{i=1}^3 L^x_4(y_i))\cap (\cap_{i=1}^3 L^y_4(y_i))\cap (\cap_{i=1}^3 L^z_4(y_i)))=2.$$
(This implies that the 13 hyperplanes $M_4(x_i), i=1,\ldots,4, L^x_4(y_i),$ $L^y_4(y_i), L^z_4(y_i), i=1,2,3$ intersect
transversally.)

\end{enumerate}
\end{lemma}

Let us prove the first statement of the lemma. Suppose that $x_1,\ldots,x_6$ are generic points in $\CP^2$. It suffices to prove that 18 linear conditions
on the space $\Pi_5$ that define the space $L(\{x_1,\ldots,x_6\})$ are independent. Suppose they are not, then $\dim_\C L(\{x_1,\ldots,x_6\}) \geq 4$. Choose
a point $x'$ such that $x_1,\ldots, x_6,x'$ are in generic position. The space $L(\{x_1,\ldots,x_6, x'\})$ would be then of dimension $\geq 1$, which is impossible,
because no curve of degree 5 can have 7 generic singular points (if the curve is irreducible, this follows from \cite[exercice 3, $\S 2$ Chapter III]{shaf}, otherwise
this is trivial).

The statement \ref{xyz2}, \ref{xyz3}, \ref{xyz5} can be proved in an analogous way. The statement \ref{xyz4} follows from B\'ezout's theorem.
$\diamondsuit$

Let us prove now that the spaces $X_i$ introduced in Proposition \ref{class} satisfy Condition \ref{three}. The case of $X_1,X_2,X_3$
is evident. If we have 4 points on a line, we can always find a quartic passing by three of these points and not passing by the fourth
point, hence the codimension (in the space of quartics) of the space of quartics passing by 4 points on a line is 4, which gives us the dimension of
$L(K), K\in X_{\ref{4ptslin}}$.

If a curve of degree 5 contains 5 singular points on a line, then this curve is defined by a polynomial of the form $f^2 g$, where $f$ is
a polynomial that defines the line, and
$g$ is a polynomial of degree 3. This gives the dimensions of all spaces $L(K), K\in X_{\ref{5ptslin}}, \ldots,X_{\ref{line}},
X_{\ref{5+1}},\ldots, X_{\ref{line+1}}, X_{\ref{5+2}},\ldots, X_{\ref{line+2}}$.

Consider a configuration $K\in X_{\ref{4nlin}}$. If no 3 of the points of $K$ are on a line, we have $\dim_\C(L(K))=21-12=9$ by the
statement \ref{xyz1} of Lemma \ref{lll}.
If $K$ contains 3 points on a line $l$, then, due to the statement \ref{xyz4} of Lemma \ref{lll}
every $f\in L(K)$ has the form $f=gh$, where $g$ is a fixed polynomial of degree 1, and $h$ is an arbitrary polynomial
of degree 4 that defines a curve that has 3 fixed intersection points with $l$ and a fixed singular point outside $l$. Using the statement
\ref{xyz5} of the same lemma (the transversality of intersection), we obtain $15-3-3=9$. The same argument
gives the dimensions of $L(K), K\in X_{\ref{4+1}}, X_{\ref{4+2}}, X_{\ref{4+3}},X_{\ref{5gen}}, X_{\ref{6gen}}$.

Note that if $l_1,l_2$ are two distinct lines in $\CP^2$,
$x_1,x_2,x_3\in l_1\setminus l_2, y_1,y_2,y_3\in l_2\setminus l_1, A\not\in l_1\cup l_2$, then the statement \ref{xyz2} of Lemma \ref{lll}
implies that
$$\dim_\C (\cap_{i=1}^3 M_3(x_i))\cap(\cap_{i=1}^3 M_3(y_i))\cap L^x_3(A)\cap L^y_3(A)\cap L^z_3(A)=1,$$
which means that the hyperplanes $M_3(x_i), M_3(y_i), i=1,2,3, L^x_3(A),$ $L^y_3(A),$ $L^z_3(A)$ intersect transversally. This
gives the dimension of $L(K), K\in X_{\ref{3+3}}, X_{\ref{3+3+int}}, X_{\ref{3+3+1}}, X_{\ref{3+3+1+int}}$. The same argument works
for $X_{\ref{6q}}, X_{\ref{6q+1}}$, except that we apply the statement \ref{xyz3} of Lemma \ref{lll} (instead of statement \ref{xyz2}).

It is easy to see that for any $K\in X_{\ref{4+3lin}},\ldots,X_{\ref{2lines}}$ the vector space $L(K)$ consists of polynomials
of the form $f^2 g^2 h$, where $f,g$ are some fixed polynomials of degree 1 that define two distinct lines, and $h$ is an arbitrary
polynomial of degree 1. Analogously, for any $K\in X_{\ref{7q}}, X_{\ref{quadric}}$ the space $L(K)$ consists of polynomials of the
form $f^2 g$, where $f$ is a fixed polynomial of degree 2 that defines a nondegenerate quadric, and $h$ is an arbitrary polynomial
of degree 1.

Consider a configuration $K\in X_{\ref{8}}$ and $f\in L(K)$. It follows from the statement \ref{xyz4} of Lemma \ref{lll} that $f=gh$, where
$g$ is a polynomial of degree 1, and $h$ is a polynomial of degree 4 that has singularities at
4 generic points $A,B,C,D$ outside the line $l$ defined by $g$. This implies that $h=h_1 h_2$, where $h_1,h_2$ are polynomials of
degree 2 that define two quadrics $Q_1, Q_2$ passing through $A,B,C,D$. $f$ must also have singularities at 4 points on $l$, hence
each of those 4 points belongs exactly to one of the quadrics $Q_1, Q_2$. It follows that $f$ is defined by $K$ up to nonzero constant.

Analogously it can be proved that for any $K\in X_{\ref{9}}, X_{\ref{10}}$ and any $f\in L(K)$, $f$ is defined by $K$ up to nonzero
constant. The cases $X_{\ref{line+3}}, X_{\ref{cp2}}$ are trivial. Thus, we have proved that the spaces $X_i$ satisfy Condition \ref{three}.

Let us prove now that these spaces satisfy Conditions \ref{six} and \ref{last}. The spaces $X_i$ satisfy
Condition \ref{xxx}. Recall that this implies (see p. \pageref{6i7}) that for the spaces $X_i$ consisting of finite configurations Condition \ref{last} follows
from Condition \ref{six}. Consider some $X_i$ that consists of finite configurations. It is immediate to check that the number of elements in all configurations
from $X_i$ is the same. Denote this number by $k$.
If we put $\eu M_k=\{(x,K)|x\in\CP^2, K\subset\CP^2, \sharp (K)=k, x\in K\}$, then $\eu M_k$ is the total space
of a fiber bundle over $B(\CP^2,k)$ (with the projection $(x,K)\mapsto K$). The triple $(\eu T_i, X_i, pr_i)$ is the restriction of this fiber bundle to $X_i$.

Thus, all spaces $X_i$ that consist of finite configurations satisfy Conditions \ref{six} and \ref{last}.

Now consider, for instance, the space $X_{\ref{2lines}}$. Note that if $G$ is a Lie group that acts smoothly on a smooth manifold $M$ and $a\in M$, there
exist submanifolds $S\subset M, S'\subset G$ such that 1). $a\in S, e\in S'$ ($e$ is the unit element of $G$),
2). $S'$ is transversal to $Stab(a)$ at $e$, 3). $S'a\subset M$ is a submanifold that intersects $S$ transversally at $a$
(here $S'x=\{gx|g\in S'\}$), 4). the map $S\times S'\to M$, defined by $(a',g)\mapsto ga', g\in S', a'\in S$ is a diffeomorphism onto an open neighborhood of $a$ in $M$. 
Put $M=X_{\ref{2lines}}, G=PGL(\CP^2)$. The action is transitive, so for any $K\in X_{\ref{2lines}}$ the above remark gives us a neighborhood $U\ni K$ and a
diffeomorphism $r:U\to S', S'\subset PGL(\CP^2)$ such that for any $K'\in U$ we have $r(K') K=K'$. Now for any $x\in K, K'\in U$, put $t(x, K')=r(K')x$. It is clear that
the map $(x,K')\mapsto (t(x,K'),K')$ is a local trivialization of $\eu T_{\ref{2lines}}$ over $U$. This trivialization satisfies Condition \ref{last}, since
all spaces $X_i$ are invariant under $PGL(\CP^2)$.

Local trivializations of tautological bundles over the spaces $X_{\ref{line}},X_{\ref{line+1}},X_{\ref{line+2}},X_{\ref{line+3}},X_{\ref{quadric}}$ can be constructed in the same way.

However, this method does not work for $X_{\ref{line+3line}}$, because the action of $PGL(\CP^2)$ on this space is no longer transitive. But we can proceed as follows:
consider $K\in X_{\ref{line+3line}}$. We have
$K=K_1\sqcup K_2$, where $K_1$ is a line, and $K_2$ consists of 3 points on another line.
Denote by $B'$ be the space of all configurations of $\CP^2$ consisting of 3 points on a line.
Let $U_1$, respectively, $U_2$ be neighborhoods of $K_1$ in $X_{\ref{line}}$, respectively, of $K_2$
in $B'$ such that the bundle $(\eu T_{\ref{line}}, X_{\ref{line}}, pr_{\ref{line}})$ is trivial over $U_1$, $(\eu T_3, X_3, pr_3)$ is trivial over
$U_2$, and for every $K_1'\in U_1', K_2'\in U_2'$ we have $K_1'\cap K_2'=\varnothing$. For $j=1,2$ let
$t_j:K_j\times U_j\to\CP^2$ be a map such that that the map $(x,K'_j)\mapsto (t_j(x,K'_j), K'_j), x\in K_j, K_j'\in U_j$ is a trivialization of the corresonding
tautological bundle over $U_j$. Put $U=\{K_1'\sqcup K_2'|K_1'\in U_1', K_2'\in U_2'\}$, and for any $K'=K_1'\sqcup K_2'\in U$ put $t(x,K')$ equal to
$t_j(x,K_j')$, if $x\in K_j, j=1,2$. It is clear that $U$ is an open neighborhood of $K$ in $X_{\ref{line+3line}}$ and that the map $(x,K')\mapsto (t(x,K'),K')$
is a trivialization of $(\eu T_{\ref{line+3line}}, X_{\ref{line+3line}}, pr_{\ref{line+3line}})$ over $U$. It follows from the construction of $t$ that for any fixed $K'\in U$,
the map $x\mapsto t(x,K')$ establishes a bijective correspondence between the
subsets of $K$ and $K'$ that belong to $\bigcup_{j\leq\ref{line+3line}}X_j$. Due to Condition \ref{xxx} it suffices to check if for maximal finite
subconfigurations of $K$ (which belong to $X_{\ref{5+3lin}},X_{\ref{6+2}},X_{\ref{7+1}}, X_{\ref{10lin}}$) and for nondiscrete subconfigurations (which belong to
$X_{\ref{line}},X_{\ref{line+1}},X_{\ref{line+2}}$).

We have shown that the spaces $X_i$ satisfy Conditions \ref{six} and \ref{last}. It remains to verify Condition \ref{five}.

Let us begin with the following lemmas.

\begin{lemma}
\label{haha}
Denote by $\Pi_d$ the vector space of all homogeneous polynomials $\C^3\to\C$ of degree $d$. The map $\Pi_d\setminus\{0\}\to 2^{\CP^2}$ that sends a polynomial
into the projectivization of the set of its zeroes is continuous.
\end{lemma}

$\diamondsuit$

\begin{corollary}
\label{c2}
The subspace of $2^{\CP^2}$ consisting of all zero sets of homogeneous polynomials of some fixed degree is closed.
\end{corollary}

{\it Remark.} In the case of real polynomials, the analogous map to the real projective plane is neither everywhere defined nor continuous
on its domain of definition.

For any $f\in\Pi_d\setminus\{0\}$ denote by $[f]$ the image of $f$ under the natural map $\Pi_d\setminus\{0\}\to (\Pi_d\setminus\{0\})/{\C^*}$.

\begin{lemma}
\label{haha1}
Suppose we have a sequence $(K_i), K_i\in 2^{\CP^2}$ and a sequence $(f_i), f_i\in\Pi_d\setminus\{0\}$, and suppose that $f_i$ has a singularity at every point of $K_i$. If
$K\in 2^{\CP^2}, f\in\Pi_d\setminus\{0\}$ are such that $\lim_{i\to\infty} K_i=K,\lim_{i\to\infty} [f]_i=[f]$, then $f$ has a singularity at every point of $K$.
\end{lemma}

$\diamondsuit$

\begin{lemma}
\label{haha2}
Suppose we have sequences $(L_i), (M_i)$ in $2^{\CP^2}$ and suppose that $K\in 2^{\CP^2}, K=\lim_{i\to\infty} (L_i\cup M_i)$. Then there exist a sequence of indices
$(i_j)$ such that $K=(\lim_{j\to\infty} L_{i_j})\cup (\lim_{j\to\infty} M_{i_j})$.
\end{lemma}

{\it Proof of Lemma \ref{haha2}.} Choose a sequence $(i_j)$ such that there exist $\lim_{i\to\infty} L_i,\lim_{i\to\infty} M_i$ and denote these limits by $L,M$ respectively.
Let $\rho$ be a metric that induces the usual topology on $\CP^2$, and let $\tilde\rho$ be the corresponding Hausdorff metric on $2^{\CP^2}$. If $A,B,C,D\in 2^{\CP^2}$, then
$\tilde\rho (A\cup B,C\cup D)\leq\tilde\rho (A,C)+\tilde\rho (B,D)$. This implies that $\tilde\rho(M_{i_j}\cup L_{i_j}, M\cup L)\leq \tilde\rho(M_{i_j}, M)+\tilde\rho
(L_{i_j},L)$. Hence $M\cup L=\lim_{j\to\infty}(M_{i_j}\cup L_{i_j})=\lim_{i\to\infty}(M_i\cup L_i)=K. \diamondsuit$

Now the verification of Condition \ref{five} becomes straitforward (in most cases). Consider, for instance, $K\in\bar X_{\ref{4+4}}$. We have $K=\lim_{i\to\infty} K_i$, since all $K_i
\in X_{\ref{4+4}}$, they can be represented as $K_i=(K_i\cap l_1^i)\cup (K_i\cap l_2^i)$, where $l_1^i,l_2^i$ are lines. Due to Lemma \ref{haha2} we can suppose that $K=(\lim_{i\to\infty}
(K_i\cap l_1^i))\cup (\lim_{i\to\infty}(K_i\cap l_2^i))$. Using Corollary \ref{c2}, we can suppose that the sequences $(l_1^i),(l_2^i)$ converge. Applying Proposition \ref{inc},
we see that $K$ is a configuration of the form ($\leq 4$ points on a line $l_1$)$\cup$($\leq 4$ points on a line $l_2$). All such configurations belong to
$\bigcup_{i=1}^{\ref{4+4}}X_i$.

The same argument works for all $X_i$ except $X_{\ref{8}},X_{\ref{9}},X_{\ref{10}}$. Let us see, what happens in these cases.
Consider, for instance, $K\in X_{\ref{8}}$. Using Lemma \ref{haha1}, we see that $K$ is included into the singular set of some polynomial $f$ of degree 5.
If this singular set is discrete, there is nothing to prove: due to Conditions \ref{first} and \ref{xxx} all subsets containing $\leq 8$ elements of discrete singular sets belong
to $\bigcup_{i=1}^{\ref{8}}X_i$. Otherwise we can do the following.

We have $K=\lim_{i\to\infty} K_i$, where all $K_i$ are of the form $(Q_1^i\cap Q_2^i)\cup (Q_1^i\cap l^i)\cup (Q_2^i\cap l^i)$, $Q^i_1,Q^i_2$ are quadrics, $l^i$ are lines. Applying
Lemma \ref{haha2} and Corollary \ref{c2}, we can suppose that $$K=(\lim_{i\to\infty}(Q_1^i\cap Q_2^i))\cup (\lim_{i\to\infty}(Q_1^i\cap l^i))\cup (\lim_{i\to\infty}(Q_2^i\cap l^i))$$
and that the sequences $(Q_1^i), (Q_2^i), l^i$ converge. Denote their limits by $Q_1, Q_2, l$ respectively.

$f$ can have the following nondiscrete singular sets: a line of multiplicity $\geq 2$, two double lines, a double line $+$ a triple line, a double nondegenerate quadric. Let us consider all
these cases.

{\it A double line.} In this case we have the following possibilities:

1). $Q_1=m_1\cup m_2, Q_2=m_1\cup m_3, m_1,m_2,m_3,l$ are 4 pairwise distinct lines. $\lim_{i\to\infty}(Q_1^i\cap Q_2^i)$ is included into a configuration of the form
(3 points on $m_1$)$\{$the point $m_2\cup m_3\}$. Thus, $K$ is included into a conficuration (the points $l\cap m_1, l\cap m_2, l\cap m_3, m_2\cap m_3$)$\cup$(3 points on $m_1$). Such a
configuration is included into a configuration from $X_{\ref{4+3}}$ or $X_{\ref{5+2}}$.

2). $Q_1=l\cup m, m\neq l, Q_2$ contains neither $l$ nor $m$. $\lim_{i\to\infty}(Q_1^i\cap l^i)$ consists of 1 or 2 points on $l$. Hence $K$ is included into a configuration of the form
$(Q_2\cap l)\cup (Q_2\cap m)\cup$(2 points on $l$), which contains $\leq 6$ points.

3). $Q_1=m\neq l, Q_2$ contains neither $l$ nor $m$. $K=(m\cap l)\cup (Q_2\cap l)\cup (Q_2\cap m)$, hence $K$ contains $\leq 5$ points.

{\it A triple line.} In this case $K$ is included into a configuration of the form (a line)$+$(a point). Hence $K$ is a subset of a configuration from $X_{\ref{8lin}}$ or $X_{\ref{7+1}}$.

{\it A line of multiplicity $\geq 4$.} $K$ is a subset of a configuration from $X_{\ref{8lin}}$.

{\it 2 double lines or a double line $+$ a triple line.} $K$ is a subset of the union of 2 lines and $\sharp (K)\leq 8$. All such configurations belong to $\bigcup_{i=1}^{\ref{4+4}}X_i$.

{\it A double nondegenerate quadric.} We have $Q_1=Q_2$ is nondegenerate. $\lim_{i\to\infty}(Q_1^i\cap Q_2^i)$ is included into a subconfiguration of the form (4 points on $Q_1$),
$\lim_{i\to\infty}(Q_1^i\cap l^i)=\lim_{i\to\infty}(Q_2^i\cap l^i)$ contains 1 or 2 points on $Q_1$. Thus, $K$ is included into a configuration from $X_{\ref{6q}}$.

We have checked Condition \ref{five} for $X_{\ref{8}}$. The spaces $X_{\ref{9}}, X_{\ref{10}}$ can be considered in a similar way.

Proposition \ref{class} is proved.
$\diamondsuit$

Now we apply Theorem \ref{conical} and Lemmas \ref{fin},\ref{fin1} to construct a conical resolution
$\sigma$ and a filtration $\varnothing\subset F_1\subset\cdots\subset F_{\ref{cp2}}=\sigma$.
The spectral sequence (\ref{spseq}) is exactly the sequence corresponding to this filtration.

Most of the columns of the sequence (\ref{spseq})
can be investigated in essentially the same way as
in the case of nonsingular quartics considered in \cite{vas1}.
We shall only discuss the columns that need a
somewhat different argument. We shall use the notations from the article
\cite{vas1} if not indicated otherwise.

\subsection{Column \ref{8}}

Let $X_{\ref{8}}$ be the space of all configurations of type \ref{8}
(see Proposition \ref{class}). From Lemma \ref{fin} we get

\begin{equation}
\label{Thom}
E^1_{\ref{8},i}=\bar{H}_{\ref{8}+i-2-7}(X_{\ref{8}}, \pm {\R}),
\end{equation}
where the local system $\pm \R$ is described in Definition \ref{def2}.
\medskip

$X_{\ref{8}}$ is naturally fibered over the space $\tilde B(\CP^2,4)$
of generic quadruples $\{A,B,C,D\} \subset \CP^2$. Let us denote by
$Y$ the fiber of this bundle, i.e. the space of all
configurations from $X_{\ref{8}}$ such that the points of intersection of
the quadrics are fixed.

\begin{lemma}
\label{X}
The term $E^2$ of the spectral sequence of the bundle
$X_{\ref{8}}\to\tilde B(\CP^2,4)$ is trivial.
\end{lemma}
The proof will take the rest of the subsection.

Denote by $L$ the space of all lines not passing through four generic points
$A, B, C, D$ in ${\CP^2}$. For any such line $l$ denote
by $Z$ the space of quadrics passing
through $A,
B, C, D$ and not tangential to $l$. The space $Z$ is homeomorphic
to ($S^2$ minus 2 points)$=\C^*$.

$Y$ is  fibered over $L$ with the fiber
$B(Z,2)=B(\C^*,2)$.

\begin{lemma}
\label{Y}
The Borel-Moore homology group of the fiber $Y$ of the bundle $X_{\ref{8}} \to
\tilde B(\CP^2,4)$ can be obtained from the spectral
sequence of the bundle $Y \to L$, whose term $E^2$ is as follows:

\begin{equation}
\label{ss2}
\begin{array}{cccc}
3&{\R^3}&{\R^3}& {\R}\\
2&\R && \\
 &2& 3 & 4
\end{array}
\end{equation}
\end{lemma}

{\it Proof.} Recall that $B(\C^*,2)$ is a fiber bundle with
the base $\C^*$ and the fiber $\C\setminus\{\pm 1\}$.
Let us study the restriction of the coefficient system
$\pm\R$ to the fiber $B(\C^*,2)$ of the bundle $Y\to L$. This system
changes its sign, when one of the points passes around zero (and the
other stands still). This corresponds to the fact that if we fix all
points in the configuration except the points of intersection of the line
and one of the quadrics, we can transpose those points. On the contrary, a loop
that transposes two quadrics, transposes two pairs of points and does not
change the sign of the coefficient system. We see that the loops of the fiber
do not change the sign of
the coefficient system, and some loop that projects onto the generator
of $\pi_1(\C^*)$ (and hence any other such loop)
does. So $\pm\R|B(\C^*,2)$ is the system $\eu A_1$ of Proposition \ref{B}.
We have $\bar{H}_2(B({\C }^*,2),\eu A_1)=\bar{H}_3(B({\C }^*,2),\eu A_1)=\R$.

The space $L$ is homeomorphic to $\C^2$ minus three generic complex lines.
We have
$\bar{H}_i(L)=\R^3$ if $i=2,3$, $\bar{H}_i(L)=\R$ if $i=4$ and $\bar{H}_i(L)=0$
otherwise. We shall complete the proof of Lemma \ref{Y} in the
following two lemmas.

\begin{lemma}
\label{loop}
Let $l(t)$ be a loop in $L$ that moves a line $l=l(0)$ around one of the points $A,B,C,D$.
Let $Z$ be the space of quadrics passing through $A,B,C,D$ and not
tangential to $l$. We can identify $Z$ with $\C^*$ (choosing an appropriate coordinate
map $z:Z\to\C^*$) in such a way that the map $Z\to Z$ induced by
$l(t)$ can be written as $z\mapsto 1/z$. If moreover
$A=(1,0), B=(-1,0), C=(0,1), D=(0,-1), l(t)=\{x=\alpha(t)\},
\alpha(t)=1+\eps e^{2\pi it}$, where $\eps=\frac{2}{\sqrt 3}-1$, then the quadrics
$q_1=\{xy=0\}$ and $q_2=\{x^2+y^2=1\}$ are preserved, and the points of
intersection of $q_1$ and $l$ are preserved, while the points of
intersection of $q_2$ and $l$ are transposed.  \end{lemma}

{\it Proof.}
Denote by $Q$ the space of quadrics passing through $A,B,C,D$.
These quadrics can be written
as follows:
$$ax^2+ay^2+bxy-a=0.$$
Such a quadric is tangential to $l(t)$ if and only if
\begin{equation}
\label{tang}
(b\alpha)^2-4a^2\alpha^2+4a^2=0.
\end{equation}

Note that if $t=0,\alpha(0)=\frac{2}{\sqrt 3}$ and this condition is
simply $a^2=b^2$.
The map $f_t:Q\to Q$ that carries the
quadrics tangential to $l$ to the quadrics tangential to $l(t)$ can be written as
$$(a,b)\mapsto(\frac{1}{2}\cdot a\sqrt{\alpha^2/(\alpha^2-1)},b).$$ If $\alpha(t)$ is as above,
$\sqrt{\alpha^2/(\alpha^2-1)}$ changes its sign, so the map $f_1(a,b)=(-a,b)$. The desired
coordinate $z\in\C^*$ is $z=(b+a)/(b-a)$. Note that $z(q_1)=1,
z(q_2)=-1, z(-a,b)=1/z(a,b)$.
The quadrics $q_1$ and $q_2$ are
preserved under any map $f_t$. The points of intersection of $q_1$ and $l$ are clearly preserved.
The statement concerning the points of intersection of $q_2$ and $l$ is verified immediately.
Lemma \ref{loop} is proved. $\diamondsuit$

Now we can describe the action of $\pi_1(L)$ on the Borel-Moore
homology of the fiber $B(Z,2)=B(\C^*,2)$ of the bundle $Y\to L$.

Lemma \ref{loop} tells us that the points of intersection of
exactly one quadric of $q_1, q_2$ are transposed. Hence the covering map of
coefficient systems $\pm\R|B(\C^*,2)\to\pm\R|B(\C^*,2)$ is minus
identity over the couple $\{\pm 1\}$. This implies that the fiber of the
coefficient system over the couple, say $\{\pm i\}$, is mapped identically.

Applying Proposition \ref{B}, we obtain immediately
that $\bar H_2(B(Z,2),\pm\R|B(Z,2))$ is multiplied by $-1$ and $\bar H_3(B(Z,2),\pm\R|B(Z,2))$ is preserved.

Thus, the 3-d line of the sequence (\ref{ss2}) contains the Borel-Moore homology of $L$ with
constant coefficients. In order to obtain the 2-nd line we must calculate the Borel-Moore
homology of $L$ with coefficients in the system $\eu L$ that changes its sign under
the action of a loop in $\CP^{2\vee}$ that embraces one of the lines corresponding to the points $A,B,C,D$.

\begin{lemma}
\label{linear}
Let $L$ be the complement of $\CP^2$ to 4 generic complex lines.
Let $f:L\to L$ be the restriction to $L$ of the projective linear map that
transposes two of those lines and preserves the other two,
and let $\tilde f:\eu L\to\eu L$ be the map that covers
$f$ and is identical over some point of $L$ that is preserved under $f$.
The Poincar{\'e} polynomial of $L$ with coefficients in $\eu L$
equals $t^2$.
The map $\tilde f$ multiplies by $-1$ the groups $H^2(L,\eu L)$ and
$\bar H_2(L,\eu L)$.
\end{lemma}

Note that since the set of fixed points of $f$ is connected if
$\tilde f$ is the identity over some fixed
point of $f$, then it is the identity over any other fixed point.

{\it Proof.}
Identify $L$ with the space
$\C^2\setminus\{z_1=0\}\cup\{z_2=0\}\cup\{z_1+z_2=1\}$.
Consider the map $p:L\to\C\setminus\{1\},p(z_1,z_2)=z_1+z_2$. Put
$A_1=p^{-1}(U_1(0)), A_2=p^{-1}(\C\setminus\{0,1\})$, where $U_1(0)$ is
the open unit disc.

The space $A_1$ is homotopically equivalent to the torus
$\{(z_1,z_2)||z_1|=\frac{1}{3}, |z_2|=\frac{1}{3}\}$. Both loops
$t\mapsto\frac{1}{3}(e^{2\pi i},1), t\mapsto\frac{1}{3}(1,e^{2\pi
i})$ act non trivially on the
fiber of $\eu L$, which implies that the cohomology groups
$H^*(A_1,\eu L)$ are zero.

The restriction
$p|A_2$ is a fibration. The restriction of $\eu L$ to the fiber
$p^{-1}(\frac{1}{2})$ in nontrivial, hence we have
$H^i(p^{-1}(\frac{1}{2}),\eu L)=\R$ if $i=1$ and is zero otherwise.
Consider the loops $\alpha:t\mapsto
1-\frac{1}{2}e^{2\pi it}, \beta:t\mapsto\frac{1}{2}e^{2\pi it}$.
It is easily checked that both of them induce the identical mapping of
$p^{-1}(\frac{1}{2})$ (hence the space $A_2$ is in fact
homeomorphic to the direct product $\C\setminus\{0,1\}\times
p^{-1}(\frac{1}{2})$). Note, moreover, that if a lifting of $\alpha$ into
$A_2$ changes the sign of $\eu L$, while a
lifting of $\beta$ does not.  Now
it is clear that the Poincar{\'e} polynomial $P(A_2,\eu L)=t^2$ and
that inclusion $A_1\cap A_2=p^{-1}(U_1(0)\setminus\{0\})\subset
A_2$ induces the isomorphism of 2-dimensional cohomology groups with
coefficients in $\eu L$.

Now consider the cohomological Mayer-Vietoris sequence corresponding to
$L=A_1\cup A_2$. Its only nontrivial terms will be

$$H^1(A_1\cap A_2,\eu L)\to H^2(L,\eu L)\to H^2(A_1,\eu L)\oplus H^2(A_2,\eu L)\to
H^2(A_1\cap A_2,\eu L)$$

The map on the right is an isomorphism, hence so is the map
on the left. So we have $P(L,\eu L)=t^2$.

The map $f$ preserves each fiber of $p$. Moreover, using the K{\"u}nneth
formula and Proposition \ref{C}, we obtain that $\tilde f$ multiplies
by $-1$ the groups $H^*(A_1\cap A_2,\eu L)$. Since the boundary
operator commutes with $\tilde f$, we obtain the statement of the
Lemma concerning the group $H^2(L,\eu L)$. The statement about the
Borel-Moore homology group follows from the Poincar{\'e} duality and the fact
that $f$ preserves the orientation.

$\diamondsuit$

Lemma \ref{Y} follows immediately from Lemma \ref{linear}. $\diamondsuit$

To complete the proof of Lemma \ref{X} we must calculate the action of
$\pi_1(\tilde B(\CP^2,4))$ on the Borel-Moore homology of $Y$ obtained
from the sequence (\ref{ss2}). This will be done in the following three
lemmas.

\begin{lemma}
\label{petlya1}
A loop $\gamma(t)$
in $\tilde B(\CP^2,4)$ that belongs to the image of $\pi_1(\tilde F(\CP^2,4))$ under the evident map $\tilde F(\CP^2,4)\to\tilde B(\CP^2,4)$
induces the identical map of the fiber $Y$ and preserves the
coefficient system $\pm\R|Y$ over it.
\end{lemma}

Note that $\pi_1(\tilde
F(\CP^2,4))\cong\Z_3$ (because $\tilde F(\CP^2,4)$ is diffeomorphic to
$PGL(\CP^2)$, which is the quotient of $SL_3(\C)$ by its center).

{\it Proof.} We can represent every $\gamma\in\pi_1(\tilde B(\CP^2,4))$ as follows:
$\gamma(t)=\{A(t),B(t),C(t),D(t)\}$, where $A(t),\ldots,D(t)$ are some paths in $\CP^2$ such that for any $t$
$A(t), B(t), C(t), D(t)$ are in generic position. If we have a $\gamma$ that
comes from $\pi_1(\tilde F(\CP^2,4))$, we have $A(0)=A(1)=A,\ldots,D(0)=D(1)=D$.
Denote by $Y_t$ the fiber
of the bundle $X_{\ref{8}}\to\tilde B(\CP^2,4)$ over $\gamma(t)$. Note that for any $t$ there exists a unique
projective linear map $M(t)$ that carries $A$ into $A(t)$, $B$ into $B(t)$, $C$ into $C(t)$ and $D$ into $D(t)$. This map induces
the map $f_t:Y_0\to Y_t$. The map $f_1$ is clearly identical. Moreover if we have a configuration $K\in Y_0$, the curve
starting at $K$ and covering $\gamma$ is $t\mapsto M(t)K$. Since $M(1)=Id_{\CP^2}$, $\gamma$ does not transpose any pair of
points from $K$. The lemma is proved.
$\diamondsuit$

\begin{lemma}
\label{petlya2}
\begin{enumerate}
\item A loop $\gamma\in\pi_1(\tilde B(\CP^2,4))$ transposing the points $A$ and $B$ induces
a bundle map $f_1:Y\to Y$. This map is covered by $\tilde f_1:\pm\R|Y\to\pm\R|Y$.

\item The corresponding map of $L$ into itself is obtained
from the projective linear map of $\CP^2$ that transposes the points $A$ and $B$ and preserves
$C$ and $D$.

\item Let $l$ be a line that is preserved under the transposition of
$A$ and $B$. The restriction of $f_1$ to the fiber over $l$ is the map
$B(\C^*,2)\to B(\C^*,2)$ induced by $z\mapsto 1/z$. The restriction of
$\tilde f_1$ to this fiber is minus identity over the couple $\{\pm i\}$.

\item The map $\tilde f_1|B(Z,2)$ multiplies by $-1$ the group $\bar H_3(B(Z,2),\pm\R|B(Z,2))$ and preserves
the group $\bar H_2(B(Z,2),\pm\R|B(Z,2))$, where $B(Z,2)$ is the fiber of the bundle $Y\to L$ over
some line from $L$ that is preserved under the transposition of $A$ and $B$.
\end{enumerate}
\end{lemma}

{\it Proof.} Proceeding as in the proof of Lemma \ref{petlya1} we obtain bundle maps
$f_t:Y_0\to Y_t$. The map $f_1:Y_0\to Y_1=Y_0$ is induced by the projective linear map preserving
$C$ and $D$ and transposing $A$ and $B$. Note that the map $\tilde
f_1:\pm\R|x\to\pm\R|x$, where $x\in Y, f_1(x)=x$ is just the map
induced by the loop $t\mapsto f_t(x)$.
Now let $A,B,C,D$ be the following
points of the affine plane $\C^2\subset\CP^2$:
$A=(1,0), B=(-1,0), C=(0,1), D=(0,-1)$. Put $l=\infty$. Then the map
$f_1$ is induced by the linear map with the matrix
$$\left(
\begin{array}{cc}
-1&0\\
0&1
\end{array}
\right)$$
This map preserves the line $l=\infty$ and transposes
the quadrics tangential to $l$. Identify $Z$ with $\C^*$, and choose such
coordinate $z\in\C^*$ that the induced map can be
written as $z\to 1/z$. Note that the quadrics in the couple that
correspond to $\{\pm i\}$ are
transposed, and the couple itself is preserved. Thus, the loop $\gamma$ transposes 3 pairs of
points over this couple, and hence the coefficient system
$\pm\R|B(Z,2)$ is multiplied by $-1$.

We have proved the first three statements. The 4-th statement
of the lemma follows immediately from Proposition \ref{B}. $\diamondsuit$

Now we can easily obtain the map of the sequence (\ref{ss2}) induced
by $\tilde f_1$.
The 3-d line of the sequence (\ref{ss2}) contains the groups $\bar H_i(L,\bar\eu H_3)$, where
$\bar\eu H_3$ is the constant local system on $L$ with the fiber
$\bar H_3(B(Z,2),\pm\R|B(Z,2))$. Since the map $\tilde f_1$
multiplies the fiber of $\bar\eu H_3$ by $-1$, it multiplies
$\bar H_4(L,\bar\eu H_3)=E^1_{4,3}$ by $-1$.  The 2-nd line of the sequence
(\ref{ss2}) contains the groups $\bar H_i(L,\bar\eu H_2)$, where $\bar\eu H_2$ is
the local nonconstant system on $L$ considered in Lemma \ref{linear}.
Due to Lemma \ref{petlya2} $\tilde f_1$ preserves the system $\bar\eu H_2$ over
some point of $L$.  We obtain from Lemma \ref{linear} that $\gamma$
multiplies $\bar H_2(L,\bar\eu H_2)=E^1_{2,2}$ by $-1$.

It is easy to see that the action of $\pi_1(\tilde B(\CP^2,4))$ on
$E^1_{2,3}$ and on $E_{3,3}^1$ of the sequence (\ref{ss2}) is
nontrivial and irreducible.

Hence the action of
$\pi_1(\tilde B(\CP^2,4))$ on $\bar H_*(Y,\pm\R|Y)$ is nontrivial
and irreducible. Recall that the universal covering space of
$\tilde B(\CP^2,4)$ is $SL_3(\C)$, and
the group $\pi_1(\tilde B(\CP^2,4))$
contains a normal subgroup
isomorphic to $\Z_3=\pi_1(\tilde F(\CP^2,4))$, the quotient being isomorphic to $S_4$.
Lemma \ref{X} follows immediately from Lemma \ref{lem3}. In fact,
putting $G=SL_3(\C), G_1=$(the subgroup of $SL_3(\C)$ generated by $e^{\frac{2}{3}\pi i}I$
($I$ is the identity matrix) and the (complexification of the) motions
of a regular tetrahedron) in Lemma \ref{lem3}
we obtain that the group $H^*(\tilde B(\CP^2,4),\eu L)=0$ if
the action of $\pi_1(\tilde B(\CP^2,4))$ on the fiber of $\eu L$ is
nontrivial and irreducible. By Poincar\'e duality $\bar H_*(\tilde B(\CP^2,4),\eu L)$ is also zero for any such $\eu L$.

\subsection{Column \ref{9}}
Recall that we denote by
$X_{\ref{9}}$ the space of configurations of type \ref{9}.
We have $$E^1_{\ref{9},i}=\bar{H}_{\ref{9}+i-2-8}(X_{\ref{9}},\pm\R).$$

$X_{\ref{9}}$ is fibered over the space $\tilde{B}(\CP^2,3)$ of all generic triples
of points $\{A,B,C\},$ $A, B, C\in \CP^2$.

\begin{lemma}
\label{24}
The spectral sequence of the bundle $X_{\ref{9}}\to\tilde B(\CP^2,3)$ looks as follows:
\begin{equation}
\label{ss7}
\begin{array}{cccccccc}
7&{\R}&&&&&&\\
6&&&{\R}&&{\R}&&\\
5&&&&&&&{\R}\\
&6&7&8&9&10&11&12\\
\end{array}
\end{equation}
and the differentials $E^2_{8,6}\to E^2_{6,7}$ and $E^2_{12,5}\to E^2_{10,6}$ are nontrivial.
\end{lemma}
The proof will take the rest of the subsection.

If we fix three lines $AB, BC,
AC$, then the intersection points of $AB$ and $BC$ with the quadric can be
chosen arbitrarily. The space of quadrics passing through these 4
points and not tangential to $AC$ is homeomorphic to ($S^2$ minus three points):
we have to exclude 2 tangential quadrics and the quadric consisting of the lines $AB$
and $BC$.

Thus, the fiber $Y$ of the bundle $X_{\ref{9}}\to\tilde B(\CP^2,3)$
is itself a fiber bundle over $B(\C^*,2)\times B(\C^*,2)$ with the fiber ($S^2$ minus
3 points). Denote the latter fiber by $Z$. The space $Z$ can be identified with
$\C^*\setminus\{1\}$.

\begin{lemma}
\label{Y'}
The spectral sequence for the Borel-Moore
homology of the bundle
$Y\to B(\C^*,2)\times B(\C^*,2)$ looks as follows:

\begin{equation}
\label{ssx}
\begin{array}{cccc}
1&\R&\R^2&\R\\
&4&5&6
\end{array}
\end{equation}
\end{lemma}

{\it Proof.} If we fix all point in a configuration from $X_{\ref{9}}$
except the points of intersection of the quadric and
the line $AC$, then we can
transpose these two points, hence the restriction of $\pm\R$ to ($S^2$ minus 3
points) is nontrivial. We have $\bar{H}_i(Z,\pm\R)=\R$ if $i=1$, and 0
otherwise.

That is why the only nontrivial line in the spectral sequence of the bundle $Y\to B(\C^*,2)\times
B(\C^*,2)$ is the first one; it contains the groups $\bar
H_*(B(\C^*,2)\times B(\C^*,2),\bar\eu H_1)$, where $\bar\eu H_1$ is the
system with the fiber $\bar H_1(Z,\pm\R|Z)$ corresponding to the
action of $\pi_1(B(\C^*,2)\times B(\C^*,2))$. The fiber of $\bar\eu H_1$
is $\R$, and, as we shall see, every element of $\pi_1(B(\C^*,2)\times
B(\C^*,2))$ multiplies the fiber of $\bar\eu H_1$ by $\pm 1$. So we can apply the K{\"u}nneth formula,
and we get
\begin{equation}
\label{product}
\bar H_*(B(\C^*,2)\times B(\C^*,2),\bar\eu H_1)=\bar
H_*(B(\C^*,2),\eu B_1)\otimes\bar H_*(B(\C^*,2),\eu B_2),
\end{equation}
where $\eu
B_1,\eu B_2$ are the restrictions of $\bar\eu H_1$ on the first and
the second factors of $B(\C^*,2)\times B(\C^*,2)$. To calculate $\eu B_1$
we fix two points in the second exemplar of $B(\C^*,2)$ and study the
action of the loops in the other exemplar of $B(\C^*,2)$ on the group $\bar H_1(Z,\pm\R|Z)$.

We put $A=(0:1:0),B=(0:0:1),C=(1:0:0)$. Put $AC=\infty$, so that the points
of the type $(z:w:1)$ belong to the affine plane
$\C^2\subset\CP^2$, and the spaces $B(\C^*,2)$ consist of couples of
nonzero points on the coordinate axes. Now fix the points
$\{(0,i),(0,-i)\}$ on the $y$-axis. Denote by $Q$ the space of quadrics
passing through $(i,0),(-i,0),(0,i),(0,-i)$. Note that the fiber
$Z$ over this quadruple is the subspace of $Q$ that consists of the
quadrics that are not tangential to $AC=\infty$ and do not equal the
union of $AB\cup AC$. Note also that the quadrics from $Q$ have the
form $$ax^2+bxy+ay^2+a.$$ Put $z=(2a-b)/2a+b)$. This identifies the space $Z\subset Q$ with $\C\setminus\{0,-1\}$.

\begin{lemma}
\label{241}
Consider the loop $\gamma(t)=(\{(\alpha(t),0),(1/\alpha(t),0)\},\{(0,i),(0,-i)\})$, where $\alpha(t)$ is a curve in
$\C\setminus \{0\}$ such that $\alpha(0)=i,
\alpha(1)=-i$. Then $\gamma$ induces the identical map of $Z$ and
preserves the points of intersection of the quadrics from $Z$ with
$\infty$.
\end{lemma}

\begin{lemma}
\label{242}
Consider the loop $\gamma(t)=(\{(ie^{\pi it},0),(-ie^{\pi
it},0)\},\{(0,i),(0,-i)\})$. The map $Z\to Z$ induced by $\gamma$ can
be written as $z\mapsto 1/z$. This map preserves the quadrics $q_1=xy$
and $q_2=x^2+y^2+1$. The points $q_1\cap\infty$ are preserved, and the
points $q_2\cap\infty$ are transposed.
\end{lemma}

The proof of these lemmas is an exercise in analytic geometry. $\diamondsuit$

Now we shall use Lemmas \ref{241} and \ref{242} to calculate the
action of $\pi_1(B(\C^*,2))$ on the Borel-Moore homology of the fiber
of the bundle $Y\to B(\C^*,2)$.

Let us note that the loop $\gamma$ considered in Lemma \ref{241} is the
loop in the fiber of the bundle $B(\C^*,2)$ over $1$.
We obtain from Lemma \ref{241} that such $\gamma$
induces the identical map of the space of $Z$, and the
points of intersection of the quadrics are preserved. So $\gamma$
transposes two points in a configuration from $Y$. Thus, the system $\pm\R|Z$
is multiplied by $-1$, and $\gamma$ multiplies by $-1$ the group
$\bar H_1(Z,\pm\R|Z)$.

Now let $\gamma$ be a loop $t\mapsto\{ie^{\pi it},-ie^{\pi it}\}$ in $B(\C^*,2)$. Note that this is exactly the
loop $a$ from Proposition \ref{B}. Applying Lemma \ref{242}, we obtain that
this loop transposes the tangential quadrics. Identify the space $Z$ of nontangential quadrics
with $\C^*$ taking the coordinate $z$ as in Lemma \ref{242}. The map $g:Z\to Z$ induced by $\gamma$ is
$z\mapsto 1/z$. Due to Lemma \ref{242} $\gamma$ transposes the points of intersection of
$q_2$ and $AC$. This implies that the map $\tilde g:\pm\R|Z\to\pm\R|Z$
induced by $\gamma$ is identity over $1=z(q_2)$ (it transposes two pairs of
points).

Now introduce another coordinate $w$ on $Z$, $w=(z-1)/(z+1)$. This identifies $Z$
with $\C\setminus\{\pm 1\}$. Since we have $w(1/z)=-w(z)$, the map $g:Z\to Z$ can be written as $w\mapsto
-w$. The map $\tilde g$ is identity over $0=w(1)$. Applying
Proposition \ref{C}, we obtain immediately that the map $\tilde g_*:\bar
H_1(Z,\pm\R|Z)\to \bar H_1(Z,\pm\R|Z)$ is minus identity.

So we see that restriction of the local system $\bar\eu H_1$ to
$B(\C^*,2)$, is in fact the system $\eu A_3$ of
Proposition \ref{B} (it changes its sign both under the action of the loops
of the fiber of $B(\C^*,2)\to\C^*$ and under the "middle
line"). Due to Proposition \ref{C} we have $\bar P(B(\C^*,2),\eu
A_3)=t^2(t+1)$. Lemma \ref {Y'} follows from formula
(\ref{product}). $\diamondsuit$

Now we shall study the action of $\pi_1(\tilde B(\CP^2,3)$ on the
group $\bar H_*(Y,\pm\R|Y)$. The fundamental group of $\tilde B(\CP^2,3)$ equals $S_3$ (since $\tilde F(\CP^2,3)$ is
simply-connected). We shall describe the map of the sequence
(\ref{ssx}) induced by the transposition of the points $A,C$. Note that there are
3 ways to represent $Y$ as a fiber bundle over $B(\C^*,2)\times B(\C^*,2)$, depending on
the choice of 2 lines, whose points of intersection with the quadric can be chosen arbitrarily.

\begin{lemma}
\begin{enumerate}
\item The transposition of the points $A$ and $C$ induces a bundle map $f_1:Y\to Y$.
The transposition of $A$ and $C$ preserves the bundle structure chosen above (we chose the lines
$AB$ and $BC$) and does not preserve the other two.

\item The corresponding map $h:B(\C^*,2)\times B(\C^*,2)\to B(\C^*,2)\times B(\C^*,2)$ is the transposition of factors.
(Recall that the first (respectively, the second) factor in this product is identified with the space of
couples of nonzero points on the $x$- (respectively, the $y$-)axis.)

\item Let $Z$ be the fiber of $Y\to B(\C^*,2)\times B(\C^*,2)$ over
$(\{(i,0),(-i,0)\},\{(0,i),(0,-i)\})$ (this point is clearly preserved
under $h$). The map $f_1:Z\to Z$ is identical, and the points of
intersection of each  quadric $q\in Z$ with the line $AC$ are transposed.
Hence a loop corresponding to the movement of any $q\in Z$ from this fiber transposes
4 pairs of points and preserves the coefficient system $\pm\R$ over
this fiber.
\end{enumerate}
\end{lemma}

{\it Proof.}
Since $\pi_1(\tilde B(\CP^2,3))=S_3$, we can take any loop $\gamma$ transposing $A$ and $C$.
Recall that $A=(1:0:0), B=(0:0:1), C=(0:1:0)\in \CP^2$. Put $A(t)=(\frac{1}{2}(1+e^{i\pi t}):
\frac{1}{2}(1-e^{i\pi t}):0), C(t)=(\frac{1}{2}(1-e^{i\pi t}):\frac{1}{2}(1+e^{i\pi t}):0)$.
Denote by $Y_t$ the fiber of the bundle $X_{\ref{9}}\to\tilde B(\CP^2,4)$ over $\gamma(t)$.
There exists a projective linear map that carries $A$ into $A(t)$, $C$ into $C(t)$ and
preserves $B$. In the affine plane $\C^2=\CP^2\setminus AC$ it looks like
$$\left(
\begin{array}{cc}
\frac{1}{2}(1+e^{i\pi t})&\frac{1}{2}(1-e^{i\pi t})\\
\frac{1}{2}(1-e^{i\pi t})&\frac{1}{2}(1+e^{i\pi t})
\end{array}
\right)$$
This map induces a map $f_t:Y_0\to Y_t$. In particular, the map $f_1$
is induced by the transposition of the axes in $\C^2$. The first and the second statements of the lemma
follow immediately. To show the 3-d statement note that the fiber $Z$ over
the point $(\{(i,0),(-i,0)\},\{(0,i),(0,-i)\})$ consists of quadrics of the type $ax^2+bxy+ay^2+a$. Such quadrics are preserved,
if we change $x$ and $y$, and their points of intersection with $AC=\infty$ are clearly transposed.
$\diamondsuit$

The map $f_1$ preserves the system $\bar\eu H_1$ corresponding
to the action of $\pi_1(B(\C^*,2)\times B(\C^*,2))$ on $\bar H_1(Z,\pm\R|Z)$, since $f_1$
induces the identical map of the fiber $Z$ over some preserved point of
$B(\C^*,2)\times B(\C^*,2)$ and preserves the restriction of the coefficient system $\pm\R$ to that fiber.

In general, suppose we have a space $A$ and a local system $\eu L$ on
$A\times A$, whose fiber is $\R$, and suppose that every element of $\pi_1(A\times
A)$ multiplies the fiber of $\eu L$ by $\pm 1$. Let $\eu L_1,\eu L_2$ be the restrictions of $\eu L$ on the first and the
second factor, and let $f:A\times A\to A\times A$ be the transposition
of factors. Suppose that $\tilde f:\eu L\to \eu L$ is a map that covers
$f$ and is identical over some point of the type $(x,x), x\in A$. Then, the
map of $\bar H_*(A\times A,\eu
L)=\bar H_*(A,\eu L_1)\otimes\bar H_*(A,\eu L_2)$ into itself can be
written as $a\otimes b\mapsto (-1)^{deg(a)deg(b)}b\otimes a$.

Applying this to our situation, we obtain that the group $\bar
H_4(B(\C^*,2)\times B(\C^*,2),\bar\eu H_1)=E^2_{4,1}$ is preserved
and $\bar H_6(B(\C^*,2)\times B(\C^*,2),\bar\eu H_1)=E^2_{6,1}$ is multiplied by $-1$.

We have clearly
$\bar{P}(\tilde B(\CP^2,3),\R)=t^{12}, \bar{P}(\tilde
B(\CP^2,3),\pm\R)=t^6$. This gives us the 5-th and the 7-th lines of ({\ref{ss7}}).

Let us calculate $\bar{P}(\tilde B(\CP^2,3),{\eu S})$, where ${\eu S}$ is the local
system corresponding to the irreducible
2-dimensional representation of $S_3$. It is easy to show that $P(\tilde
F(\CP^2,3))=(1+t+t^2)(1+t)$. Applying Corollary \ref{cor1} to the
covering $\tilde F(\CP^2,3)\to\tilde B(\CP^2,3)$, we obtain $$P(\tilde
F(\CP^2,3))=P(\tilde B(\CP^2,3),\R)+P(\tilde
B(\CP^2,3),\pm\R)+2P(\tilde B(\CP^2,3),\eu S),$$ since the regular
representation of $S_3$ contains one trivial, one alternating
representation and two exemplars of the 2-dimensional irreducible
representation. Thus, $P(\tilde B(\CP^2,3),\eu S)=t^2(1+t^2)$, and by the
Poincar{\'e} duality we obtain $\bar{P}(\tilde B(\CP^2,3),{\eu
S})=t^8(1+t^2)$.

It remains to prove that the 6-th line of the Leray sequence corresponding to the fibration
$X_{\ref 9}\to\tilde B(\CP^2,3)$ is as in ({\ref{ss7}}) (i.e., the action of $S_3$ in
$\bar{H}_5(B(\C^*,2)\times B(\C^*,2),\bar\eu H_1)=\R^2$ is irreducible), and that the differentials 
$E^2_{8,6}\to E^2_{6,7}$ and $E^2_{12,5}\to E^2_{10,6}$ are nontrivial.

To this end note that the group $SU_3$ acts almost freely on $X_{\ref 9}$ (via $SU_3\to
SU_3/\langle e^{\frac{2}{3}\pi i}I\rangle\hookrightarrow PGL(\CP^2)$, where $I$ is the identity matrix). Apply Theorem \ref{groupactions}
putting $M=X_{\ref{9}}, G=SU_3, \eu L=\pm\R|X_{\ref{9}}$. From the Leray sequence of the map
$M\to M/G$ (and from the fact that the cohomological dimension of $M/G$ is clearly finite)
we obtain that either $d_{max}-d_{min}\geq 8$ (here $d_{max},$ (respectively, $d_{min}$) is the greatest (respectively the smallest)
$i$ such that $H^i(M,\eu L)\neq 0$) or $H^*(M,\eu L)=0$.

By Poincar\'e duality, we have
either $\bar H_*(M,\eu L)=0$ or  $d'_{max}-d'_{min}\geq 8$ (here $d'_{max}$ (respectively, $d'_{min}$) is the greatest (respectively, the smallest)
$i$ such that $\bar H_i(M,\eu L)\neq 0$. Evidently if  the action of $S_3$ in
$\bar{H}_5(B(\C^*,2)\times B(\C^*,2),\bar\eu H_1)=\R^2$ is reducible or any of the differentials 
$E^2_{8,6}\to E^2_{6,7}, E^2_{12,5}\to E^2_{10,6}$ is trivial, neither statement holds. Lemma \ref{24} is proved. $\diamondsuit$

\subsection{Nondiscrete singular sets}
We are going to show that the columns of the spectral sequence
\ref{spseq} corresponding to all nondiscrete singular sets are zero.
The columns \ref{line}, and \ref{quadric} are considered in exactly the same way as in \cite{vas1}.

\begin{Prop}
\label{begemot}
Let $l$ be a line in $\CP^2$, $A_1,\ldots ,A_k, k\ge 0$ be points not on
$l$, $m$ be an integer $>1$. Then the union of all simplices in
$(\CP^2)^{*(m+k)}$ with the vertices in $A_1,\ldots ,A_k$ and $m$
vertices on $l$ has zero real homology groups modulo a point. This space is even contractible for all $k>0$.
\end{Prop}

$\diamondsuit$

Using this proposition, we can easily prove that $\partial\Lambda(K)$
has zero real homology groups modulo a point if $K\in X_i,
i=$\ref{line+1}, \ref{line+2}, \ref{line+3line}, \ref{line+3}. Let us
consider, for instance, the case $i=$\ref{line+3}. Consider a
configuration $K$ consisting of a line $l$ and 3 points $A,B,C$
outside $l$, the points $A,B,C$ do not belong to any line.

Denote the space from Proposition \ref{begemot} (that depends on
$l,m,A_1,\ldots,A_k$) by $\Lambda(l,m,A_1,\ldots,A_k)$. Note that
$\partial\Lambda(K)=L_1\cup L_2\cup L_3\cup L_4$, where
$L_1=\Lambda(l+A+B), L_2=\Lambda(l+B+C), L_3=\Lambda(l+A+C), L_4=\bigcup_{\kappa}\Lambda(\kappa)$,
where $\kappa$ runs through the set of all configurations of type
"$A,B,C+$ 4 points on $l$". Using Lemma \ref{fin1}, we conclude
that $L_4$ is homeomorphic to the space $\Lambda(l,4,A,B,C)$, which is
contractible due to Proposition \ref{begemot}.

The intersections $L_1\cap L_2, L_2\cap L_3, L_1\cap L_3$ are all
spaces of the type $\Lambda(l+\mbox{ a point not on }l)$. The
intersection $L_1\cap L_2\cap L_3$ is $\Lambda(l)$. All these spaces
are contractible. Now, the intersections $L_i\cap L_4, i=1,2,3$ are
the unions of all $\Lambda(\kappa)$, $\kappa$ running through the set
of all configurations of type "the points $(A,B)$ (reps., $(B, C)$, $(A,
C)$) outside $l$ $+$ 4 points on $l$". These spaces are
homeomorphic to the space $\Lambda(l,4,\mbox{2 points outside }l)$
and are contractible due to Proposition \ref{begemot}. Analogously, the intersections
$L_1\cap L_2\cap L_4, L_2\cap L_3\cap L_4, L_1\cap L_3\cap L_4$ are
all homeomorphic to spaces of type $\Lambda(l,4,\mbox{a point outside
}l)$ and are also contractible. Finally, the quadruple intersection
$L_1\cap L_2\cap L_3\cap L_4$ is the union of all $\Lambda(\kappa)$,
for all $\kappa=$"4 points in $l$", which is homeomorphic to $l^{*4}=\Lambda(l,4,\varnothing)$.

So we see that the spaces $L_i,i=1,\ldots,4$ have zero real homology groups modulo a point, and so do
all their intersections. This implies that real homology groups of their union
$\partial\Lambda(K)$ modulo a point are also zero.

The last case to consider is the case $K=l_1\cup l_2$, where $l_1, l_2$ a couple of
(distinct) lines (column
\ref{2lines}). Here $\partial\Lambda(K)$ is the union $L_1\cup L_2\cup L_3$,
where $L_i$ for $i=1,2$ is the union of all $\Lambda(\kappa)$, $\kappa$ running
through the set of configurations of the type "$l_i+$3 points on the
other line", and $L_3$ is the union of $\Lambda(K')$, for $K'$ running through the set $\{K'\subset l_1\cup l_2|\sharp (K')=8,
\sharp(K'\cap l_i)\geq 4, i=1,2\}$. Here $\sharp$ stands for the number of elements.

First note that the intersection $L_1\cap L_2$ is the
union of all $\Lambda(\kappa)$, $\kappa$ runs through the space of
configurations of the type "3 points on $l_1\setminus l_2$, 3 points on $l_2\setminus l_1$, the
point of intersection". It follows from Lemma \ref{fin1} that
$L_1\cap L_2$ is contractible.

The space $L_1$ admits the following filtration $\varnothing\subset\Lambda(l)\subset
M_1\subset M_2\subset M_3\subset M_4\subset M_5\subset M_6=L_1$. Here  $M_i, i=1,\ldots,5$ is
the union of all $\Lambda(\kappa),\kappa\subset K$ is a configuration of type
\ref{7+1}, \ref{line+1}, \ref{6+2}, \ref{line+2}, \ref{5+3lin}
respectively.

The space $M_1\setminus\Lambda(l_1)$ is fibered over $l_2\setminus
l_1$, the fiber over a point $A$ being homeomorphic to
$\Lambda(l_1,7,A)\setminus l_1^{*7}$. This fiber has trivial real
Borel-Moore homology. The space $M_2\setminus M_1$ is
fibered over $l_2\setminus l_1$, the fiber over a point $A$ being
homeomorphic to $\Lambda(l_1+A)\setminus\partial\Lambda(l_1+A)$, whose
Borel-Moore homology is also trivial.

The space $M_3\setminus M_2$ is fibered over the space $B(\C,2)$, the
fiber over a couple $\{A,B\}$ being homeomorphic to
$\Lambda(l_1,6,A,B)\setminus(\Lambda(l_1,6,A)\cup\Lambda(l_1,6,B))$. This
space also has trivial Borel-Moore homology.

The spaces $M_4\setminus M_3$ and $M_6\setminus M_5$ are considered in the same way as
$M_2\setminus M_1$. The space $M_5\setminus M_4$ is fibered over
$B(\C,3)$. The fiber over $\{A,B,C\}$ is homeomorphic to the space $\Lambda(l_1,5,A,B,C)\setminus(\Lambda(l_1,5,A,B)\cup
\Lambda(l_1,5,B,C)\cup\Lambda(l_1,5,A,C))$, whose Borel-Moore homology groups are zero. This implies that
$L_1$ (and also $L_2$) have zero real homology groups modulo a point.

Now consider the space $L_3$. Let $L'_3$ (respectively, $L''_3$) be the union of all $\Lambda(K')$ for $K'$ running through the set
$\{K'\subset l_1\cup l_2|\sharp (K')=8,
\sharp(K'\cap l_2)=5, \sharp(K'\cap l_1)=3\}$ (respectively, $\{K'\subset l_1\cup l_2|\sharp (K')=8,
\sharp(K'\cap l_1)=5, \sharp(K'\cap l_2)=3\}$). It is easy to see that all spaces $L'_3, L''_3, L'_3\cap L''_3,L'_3\cup L''_3$ are contractible.
The space $L_3\setminus (L'_3\cup L''_3)$ is the union of $\Lambda(K')\setminus\partial\Lambda(K')$ for $K'$ running through
the set of configurations that consist of 4 points on $l_2\setminus l_1$ and 4 points on $l_1\setminus l_2$. Using Lemma \ref{fin1}, we get
$\bar H_*(L_3\setminus (L'_3\cup L''_3),\R)=\bar H_{*-7}(B(\C,4)\times B(\C,4),\pm\R)=0$. Hence the real homology groups of $L_3$ modulo a point are zero.

We have also $L_1\cap L_3=L_3'', L_2\cap L_3=L_3', L_1\cap L_2\cap L_3=L'_3\cap L''_3$. These spaces are all contractible. This completes the proof that
real homology groups of $\partial\Lambda(K)$ modulo a point are zero, when $K$ is the union of 2 lines.


\begin{thebibliography}{99}
\bibitem{A} V. I. Arnold, {\it On some topological invariants of
algebraic functions}. Transact. (Trudy) of Moscow Math. Society, 1970, 21, p. 27-46.

\bibitem{vas1} V. A. Vassiliev, {\it How to calculate homology groups of
spaces of nonsingular algebraic projective hypersurfaces in $\CP^n$}, Proc.
Steklov Math. Inst., vol. 225, 1999, p. 121-140.

\bibitem{vas2} V. A. Vassiliev, {\it Topology of complements of
discriminants}, Phasis, Moscow, 1997 (in Russian).

\bibitem{shaf} I. R. Shafarevich, {\it Basic Algebraic Geometry 1}, Springer-Verlag, Berlin, 1994
\end{thebibliography}
\end{document}